\documentclass[11pt]{amsart}
\usepackage[english]{babel}
\usepackage[T1]{fontenc}
 
\usepackage{graphicx}
\usepackage{caption}
\usepackage{subcaption}

\usepackage[hidelinks]{hyperref}
\usepackage{indentfirst}
\usepackage{enumerate,amsmath,amssymb, mathrsfs,mathtools}
\usepackage{appendix}
\usepackage{url}
\usepackage{color}
\usepackage{accents}
\usepackage{setspace}
\usepackage{pdfpages}
\usepackage{stmaryrd}
\usepackage{amsrefs}
\usepackage{tikz}

\usepackage[margin=2.5cm]{geometry}
\allowdisplaybreaks

\renewcommand{\PrintDOI}[1]{%
  \href{https://doi.org/#1}{#1}%
}

\BibSpec{article}{%
	+{}{\PrintAuthors} {author}
	+{,}{ } {title}
	+{, }{\textit } {journal}
	+{}{ \parenthesize} {date}
		+{,  }{no. } {volume}
	+{,}{ } {pages}
	+{,}{ } {note}
	+{, }{\PrintDOI} {doi} 
}

\ExplSyntaxOn
\NewExpandableDocumentCommand{\gobblefirst}{m}
{
	\tl_tail:n { #1 }
}

\ExplSyntaxOff
\BibSpec{book}{%
    +{} {\PrintAuthors} {author}
    +{.} { }             {title}
    +{,} { }             {series}
    +{,} { vol.~}        {volume}
    +{.} { \textit}      {publisher}
    +{,} { }             {date}
    +{,} { \PrintISBNs}  {isbn}
    +{,} { \PrintDOI}    {doi}
}

\BibSpec{collection.article}{%
    +{} {\PrintAuthors}     {author}
    +{,} { }               {title}
    +{,} { \PrintConference}{conference}
    +{,} { \DashPages} {pages}
    +{} {\PrintBook}       {book}
     +{,} { }               {booktitle}
    +{,} { \PrintDateB}    {date}
    +{,} { \PrintISBNs}    {isbn}
    +{,} { \PrintDOI}      {doi}
}

\BibSpec{misc}{%
    +{} {\PrintAuthors} {author}
    +{,} { }             {title}
    +{,} { }             {organization}
    +{,} { }             {date}
    +{,} { available at \url} {url}
}

\mathcode`l="8000
\begingroup
\makeatletter
\lccode`\~=`\l
\DeclareMathSymbol{\lsb@l}{\mathalpha}{letters}{`l}
\lowercase{\gdef~{\ifnum\the\mathgroup=\m@ne \ell \else \lsb@l \fi}}%
\endgroup

\def\XXint#1#2#3{{\setbox0=\hbox{$#1{#2#3}{\int}$ }
		\vcenter{\hbox{$#2#3$ }}\kern-.6\wd0}}

\newtheorem{prop}{Proposition}
\newtheorem{thm}[prop]{Theorem}
\newtheorem{lem}[prop]{Lemma}
\newtheorem{cor}[prop]{Corollary}

\newtheorem{rem}[prop]{Remark}

\newtheorem{conj}[prop]{Conjecture}

\newcommand{\R}{\mathbb{R}}
\newcommand{\M}{\mathbb{M}}
\renewcommand{\S}{\mathbb{S}}
\renewcommand{\H}{\mathbb{H}}

\newcommand{\vol}{\mathrm{vol}}

\newcommand{\Vol}{\mathrm{Vol}}
\newcommand{\II}{\mathrm{II}}

\newcommand{\norm}[1]{\left\Vert#1\right\Vert}

\newcommand{\abs}[1]{\left\vert#1\right\vert}
\newcommand{\set}[1]{\left\{#1\right\}}
\newcommand{\brac}[1]{\left(#1\right)}
\newcommand{\scalar}[1]{\left \langle #1 \right \rangle}
\newcommand{\sscalar}[1]{\langle #1 \rangle}

\title{Alexandrov--Fenchel inequalities for weighted curvature integrals in space forms}
\author{Botong Xu}
\address{
\textnormal{Botong Xu \newline  \indent
University of Vienna \newline \indent
Faculty of Mathematics  \newline \indent
Oskar-Morgenstern-Platz 1 \newline \indent
1090 Vienna, Austria  \newline\indent 
\href{https://orcid.org/0000-0001-8878-0096}{https://orcid.org/0000-0001-8878-0096} \newline\indent	
\href{mailto:botong.xu@univie.ac.at}{botong.xu@univie.ac.at}}
}

\begin{document}
\begin{abstract}
We establish two types of Alexandrov--Fenchel inequalities for the weighted curvature integrals of static convex domains in space forms. The first type consists of a complete family of inequalities involving three weighted curvature integrals, with the second Minkowski type inequality of Xia [Calc.~Var.~Partial Differential Equations 55 (2016), Art.~96] and the weighted almost Schur lemma of Li and Xia [J.~Differential~Geom.~113 (2019), pp.~493--518] as special cases. The second type consists of a complete family of inequalities involving two weighted curvature integrals, confirming a conjectured geometric inequality by Hu and Li [Trans.~Amer.~Math.~Soc.~375~(2022), pp.~5587--5615]. The key novelty lies in observing that, via the orthogonal projection map, the weighted curvature integrals for domains in space forms admit a natural interpretation in terms of mixed volumes of Euclidean convex bodies, thus reducing the problem to the classical Alexandrov--Fenchel inequalities in convex geometry.
\end{abstract}
\maketitle
\section{Introduction}
Let $n\geq 1$, and let $\Omega$ be a bounded convex domain\footnote{Throughout this paper, bounded domains are understood to include their boundaries and are thus compact.} in $\R^{n+1}$ with $C^2$ boundary $\Sigma$. 
Let $\kappa = \brac{\kappa_1, \ldots, \kappa_n}$ denote the principal curvatures of $\Sigma$. Let $p_k: \R^n \to \R$, $k=1, \ldots, n$, denote the normalized $k$th elementary symmetric polynomials,
\begin{equation*}
p_k (a) = \binom{n}{k}^{-1} \sum_{1 \leq i_1< \cdots < i_k \leq n} a_{i_1} \cdots a_{i_k}, \quad a=\brac{a_1, \ldots, a_n}.
\end{equation*}
Set $p_0 \brac{a}=1$. With our normalization, the quermassintegrals of $\Omega$ are defined by
\begin{equation}\label{eq:quermass-smooth}
W_0 \brac{\Omega} = \brac{n+1} \Vol \brac{\Omega}, \quad 
W_k \brac{\Omega} = \int_{\Sigma} p_{k-1} \brac{\kappa} d\mu, \quad k=1, \ldots, n+1.   
\end{equation}

The Alexandrov--Fenchel inequalities for quermassintegrals in convex geometry state that 
\begin{equation}\label{AF-type1}
W_k \brac{\Omega}^2 \geq W_{k-1} \brac{\Omega} W_{k+1} \brac{\Omega}, \quad k=1, \ldots, n;
\end{equation}
see, \cite{Sch14}*{Sect.~7.4}. Equality holds in \eqref{AF-type1} if and only if $\Omega$ is a Euclidean ball. The case $k=1$ of \eqref{AF-type1} is known as the second Minkowski inequality. Another type of Alexandrov--Fenchel inequalities is given by
\begin{equation}\label{AF-type2}
W_k\brac{\Omega} \geq \brac{n+1}^{\frac{k-l}{n+1-l}}|B^{n+1}|^{\frac{k-l}{n+1-l}} W_l \brac{\Omega}^{\frac{n+1-k}{n+1-l}}, \quad 0 \leq l <k \leq n,
\end{equation}
where $B^{n+1}$ is a unit ball in $\R^{n+1}$. Equality holds in \eqref{AF-type2} if and only if $\Omega$ is a Euclidean ball. Since $W_{n+1} \brac{\Omega} = \brac{n+1} |B^{n+1}|$, \eqref{AF-type2} is a direct consequence of \eqref{AF-type1} for convex $\Omega$. However, inequalities of the form \eqref{AF-type2} hold under weaker convexity assumptions. The case $k=1$, $l=0$ reduces to the isoperimetric inequality, which holds without any convexity assumption. The case $k=2$, $l=1$ corresponds to the Minkowski inequality, for which Huisken \cite{Hui09} obtained an unpublished proof for outward-minimizing domains using inverse mean curvature flow; see also \cites{FS14} for an independent proof. For any $0 \leq l<k \leq n$, Guan and Li \cite{GL09} established \eqref{AF-type2} for domains with star-shaped, $(k-1)$-convex boundaries using the inverse curvature flow studied in \cites{Ger90, Urb90}. 

Let $\brac{\M^{n+1}, g_{\M}}$ denote either the hyperbolic space $\H^{n+1}$ or the open upper hemisphere $\S^{n+1}_+$, equipped with the standard metric. Denote by $r$ the geodesic distance from the origin $\mathcal{O}$ of $\M^{n+1}$, chosen to be the north pole when $\M^{n+1} = \S^{n+1}_+$. Define $\lambda:\R\to\R$ by
\begin{equation*}
\lambda(t)=\sinh t \quad\text{if } \M^{n+1}=\H^{n+1} \quad \text{and} \quad \lambda(t)=\sin t \quad\text{if } \M^{n+1}=\S^{n+1}_+.
\end{equation*}
Establishing analogs of the inequalities \eqref{AF-type1} and \eqref{AF-type2} in $\M^{n+1}$ has been an active area of research in recent decades. We review some closely related results in Subsection \ref{subsec:history}. 

We now introduce a generalization of the quermassintegrals in \eqref{eq:quermass-smooth} for domains in $\M^{n+1}$.
Let $\Omega$ be a bounded domain in $\M^{n+1}$ with $C^2$ boundary $\Sigma$. For any $\varepsilon \in \R$, we define the shifted principal curvatures of $\Sigma$ by $\kappa-\varepsilon := \brac{\kappa_1-\varepsilon, \ldots, \kappa_n -\varepsilon}$. Curvature flows involving shifted principal curvatures in curved ambient spaces have attracted attention in recent years; see, e.g., \cites{ACW21,BH17,HLW22,LX,WWZ23}. The constant $\varepsilon$ is typically introduced to account for the curvature of the ambient space. We introduce the following \textit{weighted shifted curvature integrals} of $\Omega$:
\begin{equation}\label{eq:WeiCurInte}
\begin{aligned}
W_0^{\lambda', \varepsilon}\brac{\Omega} &= \brac{n+1} \int_{\Omega} \lambda' \brac{r} \, d\vol, \\
W_{k}^{\lambda',\varepsilon} \brac{\Omega} &= \int_{\Sigma} \brac{\lambda'\brac{r} -\varepsilon u}p_{k-1} \brac{\kappa- \varepsilon} \, d\mu, \quad k=1, \ldots, n+1,  
\end{aligned}
\end{equation}
where $u = g_{\M} \brac{\lambda \brac{r} \partial_r, \nu}$ is the support function of $\Sigma$, and $\nu$ is the outward unit normal to $\Sigma$. When $\varepsilon =0$, we write $W_k^{\lambda'}$ for $W_k^{\lambda',0}$, a notation introduced by Hu and Li \cite{HL22} for the case $\M^{n+1} = \H^{n+1}$. The weighted shifted curvature integrals $W_k^{\lambda', \varepsilon}$ have an alternative representation:
\begin{equation*}
W_k^{\lambda',\varepsilon} \brac{\Omega} = \sum_{i=0}^k \binom{k}{i}  \brac{-\varepsilon}^{k-i} W_i^{\lambda'} \brac{\Omega}, \quad k=0,\ldots, n+1;
\end{equation*}
see Proposition \ref{prop:Wkeps-Equidef-Wk}.

A domain $\Omega$ in $\M^{n+1}$ (or its boundary $\Sigma$) is called \textit{static convex} with respect to $\mathcal{O}$ if the second fundamental form $\II_{\Sigma}$ satisfies
\begin{equation}\label{cond:StaticConvex}
\II_{\Sigma} \geq \frac{\lambda''\brac{r}}{\lambda'\brac{r}} g_{\M} \brac{\partial_r, \nu} g_{\Sigma},
\end{equation}
where $g_{\Sigma}$ denotes the induced metric on $\Sigma$. We emphasize that $\mathcal{O}$ is not assumed to lie in the interior of $\Omega$. The static convexity condition \eqref{cond:StaticConvex}, originating from the work of Brendle and Wang \cite{BW14}, arises naturally in the study of geometric inequalities involving weighted curvature integrals; see, e.g., \cites{HL22,KW23,LX19,PY23,SX19,Wu26, Xia16}. Roughly speaking, static convexity plays a role analogous to Euclidean convexity for domains in $\M^{n+1}$.

It is natural to ask whether the Alexandrov--Fenchel inequalities \eqref{AF-type1} and \eqref{AF-type2} can be generalized for the weighted shifted curvature integrals defined in \eqref{eq:WeiCurInte} of static convex domains in $\M^{n+1}$.
In the following two subsections, we review known results in this direction and present our main results in Theorem \ref{MainThm-1} and Theorem \ref{Thm:FullRangeHuLi}.

\subsection{Alexandrov--Fenchel inequalities involving three weighted curvature integrals}
By applying a generalized Reilly integral formula \cite{QX15}, Xia \cite{Xia16}*{Thm.~1.1} proved the following second Minkowski type inequality for static convex $\Omega$ in $\M^{n+1}$:
\begin{equation}
W_1^{\lambda'} \brac{\Omega}^2 \geq W_2^{\lambda'} \brac{\Omega} W_0^{\lambda'} \brac{\Omega}. \label{Ineq:Xia16}
\end{equation}
Further generalizing the integral formula, Li and Xia \cite{LX19}*{Thm.~1.9} proved a weighted almost Schur lemma for domains $\Omega \subset \M^{n+1}$ that are both convex and static convex. The result can be reformulated as
\begin{equation}
W_2^{\lambda'} \brac{\Omega}^2 \geq W_3^{\lambda'} \brac{\Omega} W_1^{\lambda'} \brac{\Omega}. \label{Ineq:LX19}
\end{equation}
Equality holds in \eqref{Ineq:Xia16} or in \eqref{Ineq:LX19} if and only if $\Omega$ is a geodesic ball in $\M^{n+1}$. Note that for \eqref{Ineq:LX19}, the convexity assumption on $\Omega$ is not explicitly stated in \cite{LX19}*{Thm.~1.9}. However, it is implicitly used in the estimate in \cite{LX19}*{Eq.~(6.2)}. A remarkable feature of \eqref{Ineq:Xia16} and \eqref{Ineq:LX19} is that, although the weighted curvature integrals depend on the position of the domain relative to the origin $\mathcal{O}$, equality holds not only for the centered ball but for all geodesic balls.

Higher-order versions of the above inequalities were conjectured in \cite{LX19}*{p.~499} and, to the best of our knowledge, remain open. Using an approach completely different from those used to prove \eqref{Ineq:Xia16} and \eqref{Ineq:LX19}, we provide an affirmative answer to this question in Theorem \ref{MainThm-1} below. The outline of its proof is presented in Subsection \ref{subsec:outline}.

\begin{thm}\label{MainThm-1}
Let $\M^{n+1}$ be either $\H^{n+1}$ or $\S^{n+1}_+$.
Let $\Omega$ be a bounded static convex domain in $\M^{n+1}$ with $C^3$ boundary $\Sigma$. Assume that $\II_{\Sigma} - \varepsilon g_{\Sigma} \geq 0$ on $\Sigma$ for some $\varepsilon \in \R$. Then the weighted shifted curvature integrals $W_k^{\lambda', \varepsilon} \brac{\Omega}$, $k=0, \ldots, n+1$, are all nonnegative and satisfy 
\begin{equation}\label{MainIneq:AlexFen}
W_k^{\lambda', \varepsilon} \brac{\Omega}^2 \geq W_{k+1}^{\lambda',\varepsilon} \brac{\Omega} W_{k-1}^{\lambda',\varepsilon} \brac{\Omega}, \quad k=1,\ldots, n.
\end{equation}
If $\II_{\Sigma} - \varepsilon g_{\Sigma}>0$ on $\Sigma$, then equality holds in \eqref{MainIneq:AlexFen} if and only if $\Omega$ is a geodesic ball in $\M^{n+1}$ that satisfies this additional condition. 
\end{thm}

\begin{rem} \label{Rmk:2} $\ $
\begin{enumerate}[(i)]
\item If $\M^{n+1} = \H^{n+1}$ and $\varepsilon \in (-\infty, 1]$, then any geodesic sphere in $\H^{n+1}$ satisfies $\II_{\Sigma} - \varepsilon g_{\Sigma} >0$. The same holds if $\M^{n+1} =\S^{n+1}_+$ and $\varepsilon \in (-\infty, 0]$.
\item \label{rmk2.2} Under the assumptions in Theorem \ref{MainThm-1}, we have
\begin{equation}\label{ineq:GeneralV-AFmain}
W_j^{\lambda',\varepsilon} \brac{\Omega}^{k-i} \geq W_i^{\lambda',\varepsilon}\brac{\Omega}^{k-j} W_{k}^{\lambda',\varepsilon}\brac{\Omega}^{j-i}, \quad 0 \leq i<j<k \leq n+1.
\end{equation}
We prove \eqref{ineq:GeneralV-AFmain} in Section \ref{sec:ProofThm}, after Remark \ref{rem:EqCase-StrictSC}.
\item When $k=1$, the inequalities \eqref{MainIneq:AlexFen} are equivalent for different choices of $\varepsilon \in \R$, as observed in \cite{LWX26}*{Sect.~4}; see also Corollary \ref{cor-LWXObser}. In particular, we can choose $\varepsilon$ such that $\II_{\Sigma} -\varepsilon g_{\Sigma}>0$ on $\Sigma$. Then Theorem \ref{MainThm-1} recovers both the inequality \eqref{Ineq:Xia16} and the characterization of its equality cases given in \cite{Xia16}. This is the main motivation for us to include the parameter $\varepsilon$ in the weighted shifted curvature integrals \eqref{eq:WeiCurInte}.
\item The additional assumption $\II_{\Sigma} -\varepsilon g_{\Sigma} >0$ for the equality case of \eqref{MainIneq:AlexFen} is technical. It allows us to apply Theorem \ref{thm:InSch14}, a result stated in \cite{Sch14}, to characterize the equality case of an Alexandrov--Fenchel inequality. We expect that this assumption can be removed. The assumption that $\Sigma$ be $C^3$ is also technical. It allows us to use the approximation lemmas developed in Subsection \ref{subsec:approx}. We expect that it can be weakened to $C^2$.
\end{enumerate}
\end{rem}

Proposition \ref{prop:ConvComp} below describes the relationship between static convexity and the curvature condition $\II_{\Sigma} - \varepsilon g_{\Sigma} \geq 0$ in Theorem \ref{MainThm-1}. Statements \eqref{prop3.3} and \eqref{prop3.4} were observed in \cite{Xia16}*{p.~3}. We prove statements \eqref{prop3.1} and \eqref{prop3.2} in Section \ref{sec:Gib}, after Proposition \ref{prop:GibH-StaticC-Convex}.
\begin{prop} \label{prop:ConvComp}
\begin{enumerate}[(i)]
\item \label{prop3.1} If $\Omega \subset \H^{n+1}$ is static convex, then $\II_{\Sigma} + g_{\Sigma} >0$.  
\item \label{prop3.2} If $\Omega \subset \H^{n+1}$ is static convex with $\mathcal{O} \in \mathrm{int} \brac{\Omega}$, then $u>0$ and $\II_{\Sigma} >0$.
\item \label{prop3.3} If $\Omega \subset \H^{n+1}$ is horo-convex, i.e., $\II_{\Sigma}- g_{\Sigma} \geq 0$, then $\Omega$ is static convex. 
\item \label{prop3.4} If $\Omega \subset \S^{n+1}_+$ is convex with $\mathcal{O} \in \mathrm{int} (\Omega)$, then $\Omega$ is static convex.
\end{enumerate}
\end{prop}

\subsection{Alexandrov--Fenchel inequalities involving two weighted curvature integrals} 
Recall the following conjecture proposed by Hu and Li in \cite{HL22}*{Conj.~1.1}.
\begin{conj}[\cite{HL22}]\label{conj:HuLi}
Let $\Omega$ be a bounded static convex domain in $\H^{n+1}$ with smooth boundary. Assume that $\mathcal{O} \in \mathrm{int} \brac{\Omega}$. Then 
\begin{equation}\label{ineq:HuLiConj}
W_k^{\lambda'} \brac{\Omega} \geq h_k \circ h_l^{-1} (W_l^{\lambda'} \brac{\Omega}), \quad 0 \leq l <k \leq n+1.
\end{equation}
Equality holds in \eqref{ineq:HuLiConj} if and only if $\Omega$ is a geodesic ball centered at the origin $\mathcal{O}$. Here, $h_k:[0, \infty) \to [0, \infty)$ is a strictly increasing function defined by $h_k \brac{t} = W_k^{\lambda'} \brac{B_t } = \abs{\S^n}\lambda \brac{t}^{n+1-k} \lambda' \brac{t}^k$, the $k$th weighted curvature integral for a geodesic ball of radius $t$, and $h_l^{-1}$ is the inverse function of $h_l$.
\end{conj}
Previous studies of the conjectured inequality \eqref{ineq:HuLiConj} have mainly relied on locally constrained curvature flows. The inequality is known to hold when
\begin{enumerate}[(i)]
\item $k=1$, $l=0$, and $\Omega$ is a bounded domain with $C^1$ boundary \cite{LX23};
\item $k=2$, $l=0$, and $\Omega$ is a bounded domain with star-shaped, strictly mean convex boundary \cite{SX19}; 
\item $1\leq k \leq n+1$, $l=0$, and $\Omega$ is static convex and contains the origin in its interior \cite{HL22}.
\end{enumerate}

We provide an affirmative answer to Conjecture \ref{conj:HuLi} in Theorem \ref{Thm:FullRangeHuLi} below.

For $\varepsilon\in (-\infty, 0]$ and $k=0,1, \ldots,n+1$, we define $h_{k, \varepsilon} : [0, \infty) \to [0, \infty)$ by $h_{k, \varepsilon} \brac{t} = W_k^{\lambda',\varepsilon} \brac{B_t}$, where $B_t \subset \H^{n+1}$ is the geodesic ball of radius $t$ centered at the origin $\mathcal{O}$. Then 
$$h_{k, \varepsilon} \brac{t} = \abs{\S^n} \lambda \brac{t}^{n+1-k} (\lambda' \brac{t} -\varepsilon \lambda \brac{t})^k,$$
which is strictly increasing in $t$.
\begin{thm}\label{Thm:FullRangeHuLi}
Let $\Omega$ be a bounded static convex domain in $\H^{n+1}$ with $C^3$ boundary $\Sigma$.  Assume that $\II_{\Sigma} - \varepsilon g_{\Sigma} \geq 0$ for some  $\varepsilon \in (-\infty,0]$. Then 
\begin{equation}\label{MainIneq:ShifHuLi-H}
W_k^{\lambda', \varepsilon} \brac{\Omega} \geq h_{k,\varepsilon} \circ h_{l, \varepsilon}^{-1} (W_l^{\lambda', \varepsilon} \brac{\Omega}), \quad 0 \leq l<k \leq n+1.
\end{equation}
Equality holds in \eqref{MainIneq:ShifHuLi-H} if and only if $\Omega$ is a geodesic ball centered at the origin $\mathcal{O}$. 
\end{thm}

\begin{rem}
Recall from Proposition~\ref{prop:ConvComp}~\eqref{prop3.2} that $\Omega \subset \H^{n+1}$ is strictly convex if it is static convex with $\mathcal{O} \in \mathrm{int} \brac{\Omega}$. Thus, the case $\varepsilon =0$ of Theorem \ref{Thm:FullRangeHuLi} confirms Conjecture \ref{conj:HuLi}. Moreover, we prove the conjectured inequality \eqref{ineq:HuLiConj} under weaker assumptions on $\Omega$.  For example, our theorem also applies to horo-convex domains $\Omega$ with $\mathcal{O} \notin \Omega$. 
\end{rem}

\subsection{Outline of the proof}\label{subsec:outline}
In the following, we present brief proofs of \eqref{MainIneq:AlexFen} in Steps \ref{Step1} to \ref{Step4} and of \eqref{MainIneq:ShifHuLi-H} in Step \ref{Step5}, both for the case $\varepsilon=0$ and for domains in $\H^{n+1}$. The proofs capture the key ideas behind the other results in this paper.

\begin{enumerate}[1.]
\item\label{Step1} The Minkowski space $\R^{n+1,1}$ is $\R^{n+2}$ equipped with the Lorentzian metric of signature $\brac{n+1,1}$. Define $\pi: \R^{n+1,1} \to \R^{n+1} \times \set{0}$ to be the orthogonal projection, where $\R^{n+1} \times \set{0}$ is naturally identified with Euclidean space $\R^{n+1}$; see Figure \ref{Fig:GibProj}. We identify $\H^{n+1}$ with the hyperboloid in $\R^{n+1,1}$ given by \eqref{ModelHyperboloid}. Let $\Omega$ be a bounded domain in $\H^{n+1}$ and set $K_{\Omega} = \pi \brac{\Omega}$.  A key property of the projection $\pi$ is that $W_0^{\lambda'} \brac{\Omega} = \brac{n+1} \Vol \brac{K_{\Omega}}$, where $\Vol \brac{\,\cdot \,}$ denotes the Euclidean volume. This identity was used by the author and Li \cite{LX23} to establish \eqref{ineq:HuLiConj} for all bounded domains in $\H^{n+1}$ when $k=1$, $l=0$. 

\item \label{Step2}Let $\Omega_t = \set{X \in \H^{n+1};~ \mathrm{dist}_{g_{\H}} \brac{X, \Omega} \leq t}$ denote the outer parallel set of $\Omega$ at distance $t \geq 0$. 
We derive a Steiner type formula for $W_0^{\lambda'} \brac{\Omega_t}$  in Proposition \ref{prop:SteinerFormulaWeightedCurvature}, from which we see that $W_k^{\lambda'} \brac{\Omega}$ is, up to a constant factor, the coefficient of the term $\lambda' \brac{t}^{n-k+1} \lambda \brac{t}^k$. 

\item\label{Step3} Let $K_{\Omega_t} = \pi \brac{\Omega_t}$. In light of the above two steps, we proceed to investigate formulas for $\Vol \brac{K_{\Omega_t}}$. The boundary $\Sigma_t = \partial \Omega_t$ is given by $\Sigma_t =\{\lambda'\brac{t} X+ \lambda \brac{t} \nu \brac{X};~ X \in \Sigma\}$. Let $L_{\Omega}$ denote the domain enclosed by  $\pi \brac{\nu \brac{\Sigma}}$ when $\pi \brac{\nu \brac{\Sigma}}$ is a closed embedded hypersurface. One may therefore expect that $K_{\Omega}$ and $L_{\Omega}$ are convex bodies in $\R^{n+1}$ and that $K_{\Omega_t}$ is their Minkowski linear combination. Precisely,
\begin{equation*}
K_{\Omega_t} = \lambda'\brac{t} K_{\Omega} + \lambda \brac{t} L_{\Omega}.
\end{equation*}
We prove in Lemma \ref{lem:DecompKt-H} that this is indeed the case when $\Omega$ is both convex and static convex. By a theorem of Minkowski, $\Vol \brac{K_{\Omega_t}}$ is a polynomial in $\lambda\brac{t}$ and $\lambda'\brac{t}$.

\item \label{Step4} Recall from Step \ref{Step1} that $(n+1)\Vol \brac{K_{\Omega_t}} = W_0^{\lambda'} \brac{\Omega_t}$. By comparing the coefficients of the polynomial expansion of $\Vol\brac{K_{\Omega_t}}$ in Step \ref{Step3} with the Steiner type formula for $W_0^{\lambda'} \brac{\Omega_t}$ mentioned in Step \ref{Step2}, we find that
\begin{equation}\label{KeyRelation}
W_{k}^{\lambda'} \brac{\Omega} = \brac{n+1} V \brac{K_{\Omega}, n+1-k; L_{\Omega}, k}, \quad  k=0, \ldots, n+1.
\end{equation}
Here, $V \brac{K_{\Omega}, n+1-k; L_{\Omega}, k}$ denotes the mixed volume of $n+1-k$ copies of $K_{\Omega}$ and $k$ copies of $L_{\Omega}$; see \eqref{fml:ClassicalSteiner} for the definition. Hence, the case $\varepsilon=0$ of inequality \eqref{MainIneq:AlexFen} follows from the Alexandrov--Fenchel inequality \eqref{ineq:AFReal} for mixed volumes in Euclidean convex geometry. By Theorem \ref{thm:InSch14}, a result stated in \cite{Sch14}, the equality case of \eqref{MainIneq:AlexFen} requires $K_{\Omega}$ and $L_{\Omega}$ to be \textit{homothetic}, which holds if and only if $\Omega$ is a \textit{geodesic ball} in $\H^{n+1}$.

\item \label{Step5} Unlike the Euclidean case, the highest-order weighted curvature integral $W_{n+1}^{\lambda'}\brac{\,\cdot \,}$ depends on the domain in $\H^{n+1}$. Hence, the inequality \eqref{ineq:HuLiConj} does not directly follow from \eqref{MainIneq:AlexFen}.
To establish \eqref{ineq:HuLiConj}, we first prove the case $k=n+1$, $l=n$ using the projection $\pi$. Then we use \eqref{MainIneq:AlexFen} and an induction argument to show that \eqref{ineq:HuLiConj} holds for the cases $k=n+1$, $0\leq l\leq n$. Finally, the case $\varepsilon=0$ of Theorem \ref{Thm:FullRangeHuLi} follows by combining these cases of \eqref{ineq:HuLiConj} with \eqref{MainIneq:AlexFen}.
\end{enumerate}

\subsection{Related geometric inequalities and main contributions} \label{subsec:history}
By applying his integral formula to a solution of a suitable elliptic equation, Reilly \cite{Rei80} established a generalization of the second Minkowski inequality for bounded convex domains in manifolds with nonnegative Ricci curvature. Using generalized Reilly integral formulas, Xia \cite{Xia16} established \eqref{Ineq:Xia16}, while Li and Xia \cite{LX19} established \eqref{Ineq:LX19}. To the best of our knowledge, these are the only known sharp Alexandrov--Fenchel type inequalities in space forms of the form \eqref{AF-type1}. Motivated by an idea from \cite{Xia16}, a related inequality involving linear combinations of weighted curvature integrals in $\H^{n+1}$ was established in \cite{LX}*{Sect.~6}.

Many Alexandrov--Fenchel type inequalities of the form \eqref{AF-type2} have been established  in space forms using curvature flows. The general strategy is to construct a curvature flow along which either a suitable combination of geometric quantities of the evolving hypersurface is monotone, or one quantity is preserved while another is monotone. This reduces the proof to analyzing the limiting behavior of the flow. Many of these inequalities involve $W_k^{\lambda',\varepsilon} \brac{\Omega}$. For $\varepsilon=0$, these include results obtained using inverse mean curvature flow \cites{BHW16, dG16}, locally constrained curvature flows \cites{DL25, HL22, PS25}, and the conformal flow introduced in \cite{Bre13} and utilized in \cite{GWW15}. For $\varepsilon \neq 0$, related inequalities involving shifted principal curvatures have also been established; see, e.g., \cites{HLW22,HWZ24,LWX26}. Further inequalities with weights more general than $\lambda'\brac{r}$ have also been established; see, e.g., \cites{KM15, KW23, KW26, PL25, SW26, WZ23, Wu26}.

Previous approaches to Alexandrov--Fenchel type inequalities involving weighted curvature integrals in space forms and more general manifolds have mainly relied on analytic tools, such as Reilly integral formulas and curvature flows. These approaches yield geometric inequalities that capture the influence of ambient curvature, sometimes under weaker convexity assumptions than those required in classical convex geometry. By contrast, we focus on understanding the geometric meaning of these weighted curvature integrals. The key ingredient is the identity \eqref{KeyRelation}, which connects these integrals to the mixed volumes of Euclidean convex bodies via the projection map $\pi$. This identity allows us to reduce the previously open cases addressed by Theorem \ref{MainThm-1} to the Alexandrov--Fenchel inequalities in Euclidean space. Moreover, the identity plays a crucial role in the proof of Theorem \ref{Thm:FullRangeHuLi}.

The rest of this paper is organized as follows. In Section \ref{sec:Prel}, we review concepts and results in convex geometry and study the properties of hypersurfaces in space forms. In Section \ref{sec:Gib}, we investigate the relationship between hypersurfaces in space forms and their images in Euclidean space via the orthogonal projection $\pi$. In Section \ref{sec:Steiner}, we establish a Steiner type formula for the weighted volume. In Section \ref{sec:Decomp}, we study the decomposition of the projection of the outer parallel set of a domain. In Section \ref{sec:ProofThm}, we prove Theorems \ref{MainThm-1} and \ref{Thm:FullRangeHuLi}. In Appendix \ref{App:A}, we explain why Theorem \ref{Thm:FullRangeHuLi} has no analog in $\S^{n+1}$. In Appendix \ref{App:B}, we derive an alternative form of the Steiner formulas in space forms.

\subsection*{On the use of artificial intelligence (A.I.)} Before this work, during a discussion with GPT-5.6 Sol, the AI suggested a Steiner type formula for $W_0^{\lambda'}$ along with a proof. This was subsequently formulated as the case $\M^{n+1} = \H^{n+1}$, $\varepsilon =0$ of Proposition \ref{prop:SteinerFormulaWeightedCurvature} and independently verified by the author. The author has not been able to find this form of the Steiner formula in the literature. If this form of the formula is already known, the author would be grateful to learn of a reference. The asymptotic behavior of the weighted volume $W_0^{\lambda'}$ along outer parallel sets was discussed in \cite{Xia16}*{Remark~2.1}. A Steiner type formula was obtained in a different but equivalent form for horo-convex domains in the author's joint work with Li \cite{LX}*{Thm.~6.4}. 

This Steiner type formula reminded the author of his earlier joint works \cites{LX23, HLWX25} on orthogonal projection and served as the starting point for this work. The main results, their proofs, and the entire development of the paper were completed by the author. After the original draft was completed, GPT-6 Astra was used to review and refine the manuscript. The author takes full responsibility for the paper.

\subsection*{Acknowledgments} 
The author thanks Michael Eichmair and Thomas Koerber for reading an earlier version of the manuscript and offering helpful suggestions that improved its readability.

Part of this work was conducted while the author was a postdoctoral researcher at Technion-I.I.T.
The research leading to these results is part of a project that has received funding from the European Research Council (ERC) under the European Union's Horizon 2020 research and innovation programme (grant agreement No 101001677). This research was funded in part by the Austrian Science Fund (FWF) [10.55776/PAT9828924].

\section{Preliminaries}\label{sec:Prel}
\subsection{Convex bodies in Euclidean space}\label{subsec:2.1}
In this subsection, we collect concepts and results from convex geometry; see Schneider's monograph \cite{Sch14}. 

A convex body in $\R^{n+1}$ is a nonempty, compact convex set. Here, we do not require convex bodies to have nonempty interior. 
Given a convex body $K$, its support function $h_K: \S^n \to \R$ is defined by
\begin{equation*}
h_K \brac{\theta} = \sup_{x \in K} \sscalar{x, \theta}, \quad \theta \in \S^n.
\end{equation*}
The set of convex bodies is closed under Hausdorff convergence. Moreover, convergence in the Hausdorff metric is equivalent to uniform convergence of the corresponding support functions on $\S^n$.

Denote by $g_{\S^n}$ and $D$ the canonical metric and the Levi--Civita connection on the unit sphere $\S^n \subset \R^{n+1}$, respectively. 

If $\partial K$ is $C^2$ and strictly convex in the sense that its second fundamental form is positive definite on $\partial K$, then $h_K \in C^2 \brac{\S^n}$ and $D^2 h_K + h_K g_{\S^n}$ is positive definite on $\S^n$. Moreover, for each $\theta \in \S^n$, the function $\sscalar{x, \theta}$ attains its maximum over $K$ at the unique point $x_{\theta} = D h_K \brac{\theta}+ h_K \brac{\theta} \theta$ on $\partial K$. It follows that $\partial K = \set{Dh_K \brac{\theta}+ h_K \brac{\theta} \theta;~ \theta \in \S^n}$. Conversely, given $f \in C^2 \brac{\S^n}$ with $D^2 f + f g_{\S^n}$ positive definite on $\S^n$, there exists a unique convex body $L$ whose support function is $f$. The boundary $\partial L$ is $C^2$ and strictly convex.

Given $m$ convex bodies $K_1, \ldots, K_m$ in $\R^{n+1}$ and $m$ nonnegative coefficients $a_1, \ldots, a_m$, their Minkowski combination is defined by 
\begin{equation*}
\sum_{i=1}^m a_i K_i = \set{\sum_{i=1}^m a_i x_i;~ x_i \in K_i, \, i=1, \ldots, m},
\end{equation*}
which is a convex body with support function $\sum_{i=1}^m a_i h_{K_i}$. 
A theorem of Minkowski states that there exists a unique nonnegative symmetric $(n+1)$-linear form $V$ with respect to the Minkowski sum such that $V \brac{K, \ldots, K} = \Vol \brac{K}$ for any convex body $K$ in $\R^{n+1}$. Then
\begin{align}
\Vol \brac{\sum_{i=1}^m a_i K_i } &= V \brac{\sum_{i_1=1}^m a_{i_1} K_{i_1}, \ldots, \sum_{i_{n+1} =1}^m a_{i_{n+1}} K_{i_{n+1}} } \nonumber \\
&= \sum_{i_1, \ldots, i_{n+1} =1}^m V \brac{K_{i_1}, \ldots, K_{i_{n+1}}} a_{i_1} \cdots a_{i_{n+1}}, \label{eq:MixedVolume}
\end{align}
where $V \brac{K_{i_1}, \ldots, K_{i_{n+1}}}$ is called the mixed volume of $K_{i_1}, \ldots,  K_{i_{n+1}}$.

The following theorem characterizes the positivity of mixed volumes; see \cite{Sch14}*{Thm.~5.1.8}.
\begin{thm}[\cite{Sch14}]\label{thm:PosMixVol}
Let $K_1, \ldots, K_{n+1}$ be convex bodies in $\R^{n+1}$. Then $V \brac{K_1, \ldots, K_{n+1}}>0$ if and only if there are segments $S_i \subset K_i$, $i=1, \ldots, n+1$, with linearly independent directions.
\end{thm}

Taking $m=2$, $(a_1, a_2) = \brac{1, t}$ and $\brac{K_1, K_2} = \brac{K,L}$ in \eqref{eq:MixedVolume}, we have
\begin{equation}
\Vol \brac{K+ tL} = \sum_{i=0}^{n+1} \binom{n+1}{i} V \brac{K, n+1-i; L, i} t^i, \quad t \geq 0,\label{fml:ClassicalSteiner}
\end{equation}
where $V \brac{K, n+1-i; L, i}$ denotes $V \brac{K, \ldots, K, L, \ldots, L}$ with $n+1-i$ copies of $K$ and $i$ copies of $L$. 

Differentiating both sides of \eqref{fml:ClassicalSteiner} at $t=0$ yields
\begin{equation}\label{eq:MixVolume-2Body}
\left. \frac{d}{dt} \right|_{t=0^+} \Vol \brac{K+tL} = \brac{n+1} V \brac{K, n;L,1}.
\end{equation}
There exists a unique finite Borel measure $dS_K$ on $\S^n$ called the surface area measure of $K$ such that the right-hand side of \eqref{eq:MixVolume-2Body} can be expressed as
\begin{equation}\label{eq:AreaMeasure}
\brac{n+1}V \brac{K, n; L, 1} =  \int_{\S^n} h_L \, dS_K.
\end{equation}
Note that $\mathrm{Area} \brac{\partial K} = \int_{\S^n} dS_K$ if $K$ has nonempty interior. Taking $L =K$ in \eqref{eq:AreaMeasure}, we have
\begin{equation} \label{eq:VolMeasure}
\brac{n+1}\Vol \brac{K} = \int_{\S^n} h_K \, dS_K.
\end{equation}
The measure $dV_K :=h_K dS_K/ \brac{n+1}$ is known as the cone-volume measure of $K$ when $\mathcal{O} \in \mathrm{int} \brac{K}$. 
 
Let $K$, $L$, $C_1, \ldots, C_{n-1}$ be convex bodies in $\R^{n+1}$. The Alexandrov--Fenchel inequality states that
\begin{align}
V\brac{ K, L, C_1, \ldots, C_{n-1} }^2 \geq V\brac{K, K, C_1, \ldots, C_{n-1} } V \brac{ L, L, C_1, \ldots, C_{n-1} };  \label{Alexandrov-Fenchel Inequality}
\end{align}
see \cite{SvH19} for an elegant proof. For $k=1, \ldots, n$, taking the convex bodies $C_1, \ldots, C_{n-1}$ in \eqref{Alexandrov-Fenchel Inequality} as $(n-k)$ copies of $K$ and $\brac{k-1}$ copies of $L$, we obtain
\begin{equation}\label{ineq:AFReal}
V \brac{K, n+1-k; L, k}^2 \geq V \brac{K, n+2-k;L, k-1} V \brac{K, n-k; L, k+1}.
\end{equation}
If $L$ is a unit ball $B^{n+1}$ in $\R^{n+1}$, then \eqref{ineq:AFReal} becomes
\begin{equation}\label{ineq:Type1inSectPre}
W_k \brac{K}^2 \geq W_{k-1}\brac{K} W_{k+1} \brac{K}, \quad k=1, \ldots, n,
\end{equation}
where $W_k\brac{K}$ is called the $k$th quermassintegral of $K$, defined by
\begin{equation*}
W_k \brac{K} = \brac{n+1} V \brac{K, n+1-k;B^{n+1},k}, \quad k=0, \ldots, n+1.
\end{equation*}
In particular, $W_{n+1} \brac{K} = \brac{n+1} |B^{n+1}|$ for any convex body $K$. Combining this fact with \eqref{ineq:Type1inSectPre}, we obtain \eqref{AF-type2} for convex bodies.
Note that the quermassintegrals can be expressed by \eqref{eq:quermass-smooth} when the boundary $\partial K$ is $C^2$.

Note that our indexing and normalization conventions for the quermassintegrals differ from those used in some parts of the literature. For instance, $W_k \brac{K}$ may be defined as $V (K,k;B^{n+1},n+1-k)$ or $V (K,n+1-k;B^{n+1},k)$. Our convention is consistent with the definition of $W_k^{\lambda'} \brac{\Omega}$ in \eqref{eq:WeiCurInte}.

The characterization of the equality cases of the Alexandrov--Fenchel inequality \eqref{Alexandrov-Fenchel Inequality} is a deep problem in convex geometry that remains open in full generality; see, e.g., \cite{SvH23}*{Sect.~16.1}. However, it follows from Theorem \ref{thm:InSch14} (see \cite{Sch14}*{Theorem~7.6.8}) below that for convex bodies $K$ and $L$ with $C^2$ boundaries, equality in \eqref{ineq:AFReal} implies that $K$ and $L$ are homothetic.

\begin{thm}[\cite{Sch14}]\label{thm:InSch14}
Let $K,L, C_1, \ldots,C_{n-1}$ be convex bodies in $\R^{n+1}$ with nonempty interiors. Assume that $C_1, \ldots, C_{n-1}$ are smooth.
If equality holds in \eqref{Alexandrov-Fenchel Inequality}, namely,
\begin{equation*}
V\brac{ K, L, C_1, \ldots, C_{n-1} }^2 = V\brac{K, K, C_1, \ldots, C_{n-1} } V \brac{ L, L, C_1, \ldots, C_{n-1} },
\end{equation*}
then $K$ and $L$ are homothetic, i.e., there exist $\alpha >0$ and $v \in \R^{n+1}$ such that $K= \alpha L+ v$. Here, a convex body is called smooth if it has a unique supporting hyperplane at each boundary point. In particular, every convex body with nonempty interior and $C^1$ boundary is smooth in this sense.
\end{thm}

\subsection{Space forms}
A space form is a complete, connected Riemannian manifold of constant sectional curvature. Up to scaling, every space form is a quotient of hyperbolic space, Euclidean space, or the sphere by a subgroup of isometries.
Denote by $\mathbb{M}^{n+1}(-1) = \mathbb{H}^{n+1}$, $\mathbb{M}^{n+1}(0) =\mathbb{R}^{n+1}$, and $\mathbb{M}^{n+1}(1) =\mathbb{S}^{n+1}_+$ the $(n+1)$-dimensional hyperbolic space, Euclidean space, and upper hemisphere, respectively. For $c \in \set{-1,0,1}$, the metric $g_{\M}$ of $\M^{n+1}(c)$ in geodesic polar coordinates is given by
\begin{equation*}
g_{\M} = dr^2+ \lambda\brac{r}^2 g_{\S^n},
\end{equation*}
where
\begin{equation*}
\lambda(r) = \left\{
\begin{aligned}
&\sinh r, \quad &r \in [0, \infty), \quad &{\rm when} \ c=-1, \\
&r, \quad &r \in [0, \infty), \quad &{\rm when} \ c=0, \\
&\sin r, \quad &r \in [0, \pi/2), \quad &{\rm when} \ c=1,
\end{aligned}
\right.
\end{equation*}
and $r$ denotes the geodesic distance from a fixed point $\mathcal{O} \in \M^{n+1} (c)$, chosen to be the north pole when $c=1$. Note that $\lambda'' + c \lambda =0$ and $\brac{\lambda'}^2+c \lambda^2=1$.

Next, we recall the hyperboloid model of hyperbolic space. The Minkowski space $\R^{n+1,1}$ is an $(n+2)$-dimensional vector space equipped with a Lorentzian metric
\begin{equation*}
\sscalar{X, Y}_1= \sum_{i=1}^{n+1} x_i y_i - x_{n+2} y_{n+2}.
\end{equation*}
A vector $X \in \R^{n+1,1}$ is called timelike, lightlike, or spacelike if $\sscalar{X,X}_1<0$, $\sscalar{X,X}_1=0$, or $\sscalar{X,X}_1>0$, respectively.
For convenience, we write $X = \brac{x, x_{n+2}}$, where $x \in \R^{n+1}$, and we identify $\R^{n+1}$ with the hyperplane $\R^{n+1} \times  \set{0}$ in $\R^{n+1,1}$. The hyperboloid model of hyperbolic space $\H^{n+1}$ is given by
\begin{equation}\label{ModelHyperboloid}
\H^{n+1} = \set{X= \brac{x,x_{n+2}} \in \R^{n+1,1};~\sscalar{X,X}_1 = -1, \, x_{n+2} >0 }.
\end{equation}
Setting the origin $\mathcal{O}= \brac{0,1} \in \R^{n+1} \times \R$, a point with polar coordinates $\brac{r, \theta}$ in $\H^{n+1}$ is identified with $\brac{\lambda\brac{r} \theta, \lambda' \brac{r}}$ on the hyperboloid. It follows that $\lambda'\brac{r (X)} = x_{n+2}= - \sscalar{X, (0,1)}_1$.
For any $X \in \H^{n+1}$ and unit vector $v \in T_X \H^{n+1}$, the geodesic $\gamma \brac{t}$ satisfying $\gamma \brac{0} = X$ and $\gamma' \brac{0} = v$ is given by
\begin{equation}\label{eq:geoH}
\gamma \brac{t} = \lambda' \brac{t} X + \lambda \brac{t} v.
\end{equation}
Thus, for any two points $X_1, X_2 \in \H^{n+1}$ at geodesic distance $d$, we have
\begin{equation}\label{fml:GeoDis-H}
\lambda' \brac{d} = - \sscalar{X_1, X_2}_1.
\end{equation}

We now recall the corresponding realization of the upper hemisphere $\S^{n+1}_+$ given by 
\begin{equation}\label{ModelHemisphere}
\S^{n+1}_+ =\set{P= \brac{p, p_{n+2}} \in \R^{n+2};~ \sscalar{P,P} =1, \, p_{n+2} >0},
\end{equation}
where we use $\sscalar{\, \cdot \,, \, \cdot \,}$ to denote the Euclidean inner product.
Setting the origin $\mathcal{O} = \brac{0,1} \in \R^{n+1} \times \R$, a point with polar coordinates $\brac{r, \theta}$ in $\S^{n+1}_+$ is identified with $\brac{\lambda \brac{r} \theta, \lambda' \brac{r}}$ on the upper hemisphere. Thus $\lambda'\brac{r \brac{P}} = \sscalar{P, (0,1)}$.
For any $P \in \S^{n+1}_+$ and unit vector $v \in T_P \S^{n+1}$, the geodesic $\gamma \brac{t}$ satisfying $\gamma \brac{0} = P$ and $\gamma' \brac{0} = v$ is given by
\begin{equation}\label{eq:geoS}
\gamma \brac{t} = \lambda' \brac{t} P + \lambda \brac{t} v.
\end{equation}
For any points $P_1, P_2 \in \S^{n+1}_+$ at geodesic distance $d$, we have
\begin{equation}\label{fml:GeoDis-S}
\lambda' \brac{d} = \sscalar{P_1, P_2}.
\end{equation}

\subsection{Hypersurfaces in space forms}
In the following, $\M^{n+1}$ denotes $\H^{n+1}$, $\R^{n+1}$, or $\S^{n+1}_+$, and we write its metric $g_{\M}$ as $\sscalar{\, \cdot \, , \, \cdot \,}$.

Let $\Omega$ be a bounded domain in $\M^{n+1}$ with $C^2$ boundary $\Sigma$. We denote by $g_{\Sigma}$ and $\nabla$ the induced metric and the Levi-Civita connection on $\Sigma$, respectively. Denote by $\nu_{\Sigma}$ the outward unit normal to $\Sigma$. The second fundamental form $\II_{\Sigma}$ of $\Sigma$ is defined by $\sscalar{dX, d\nu_{\Sigma}}$. In a local coordinate system $(y^1, \ldots, y^n)$ of $\Sigma$, we denote $\brac{g_{\Sigma}}_{ij}= g_{\Sigma}\brac{\partial_i, \partial_j}$ and $\brac{\II_{\Sigma}}_{ij} = \II_{\Sigma} \brac{\partial_i, \partial_j}$, $i,j =1, \ldots,n$. The corresponding Weingarten map is $ g_{\Sigma}^{-1} \II_{\Sigma}$, and its eigenvalues $\kappa = \brac{\kappa_1, \ldots, \kappa_n}$ are called the principal curvatures of $\Sigma$. The principal curvatures of a geodesic sphere in $\M^{n+1}$ of radius $\rho$ are all equal to $\lambda'\brac{\rho}/ \lambda \brac{\rho}$. For any $\varepsilon \in \R$, we write $\kappa - \varepsilon = \brac{\kappa_1-\varepsilon, \ldots, \kappa_n -\varepsilon}$ and call $\kappa -\varepsilon$ the shifted principal curvatures of $\Sigma \subset \M^{n+1}$. Define the support function $u_{\Sigma}: \Sigma \to \R$ of $\Sigma$ by $u_{\Sigma} \brac{x} = \sscalar{\lambda \brac{r} \partial_r, \nu_{\Sigma}}$. We write $u$ instead of $u_{\Sigma}$ if there is no chance of confusion.

Let $\Omega$ be a bounded domain in $\R^{n+1}$ with strictly convex boundary $\Sigma$, in the sense that $\II_{\Sigma}>0$. Let  $G: \Sigma \to \S^n$ denote the Gauss map with $G \brac{x}= \nu_{\Sigma} \brac{x}$. Then $u_{\Sigma} \brac{x} = h_{\Omega} \brac{\nu_{\Sigma} \brac{x}}$, equivalently, 
\begin{equation}\label{eq:Rel-Supports}
h_{\Omega} \brac{\theta} = u_{\Sigma} \brac{G^{-1} \brac{\theta}}, \quad \theta \in \S^n.
\end{equation}
For each $\theta \in \S^n$, we can identify $T_{\theta} \S^n$ with $T_{G^{-1} \brac{\theta}} \Sigma$ since they are both perpendicular to $\theta$. Let $\set{e_1, \ldots, e_n}$ be a local orthonormal frame on $\S^n$. Using this identification, we have
\begin{equation}\label{fml:2FF-D2h}
\brac{\II_{\Sigma} \brac{G^{-1} (\theta) }}_{ik} \brac{D_k D_j h_{\Omega} \brac{\theta}+ h_{\Omega} \brac{\theta} \delta_{kj}} = \delta_{ij}, \quad \theta \in \S^n;
\end{equation}
see, e.g., \cite{Guan18}*{Eqs.~(1.3)--(1.4)}.
Throughout this paper, we use the Einstein summation convention.

It is a standard fact that $\lambda\brac{r} \partial_r$ is a conformal Killing vector field in $\M^{n+1}$. For any two vector fields $U$ and $V$ in $\M^{n+1}$, it holds
\begin{equation}\label{eq:ConfKilling}
\sscalar{\nabla^{\M}_U \brac{\lambda \brac{r} \partial_r}, V} = \lambda' \brac{r} \sscalar{U, V}.
\end{equation}
Using the divergence theorem, we have 
\begin{equation}\label{eq:diverW0}
\int_{\Sigma} u \, d\mu = \int_{\Sigma} \sscalar{\lambda \brac{r} \partial_r, \nu} \, d\mu = \brac{n+1} \int_{\Omega} \lambda' \brac{r} d\vol= W_0^{\lambda'} \brac{\Omega}.
\end{equation}
On the other hand, using $\lambda''+c\lambda =0$ and \eqref{eq:ConfKilling}, we have
\begin{equation}
\nabla_j \nabla_i \lambda'\brac{r} = - c \nabla_j \sscalar{ \lambda \brac{r} \partial_r, e_i }= -c \brac{\lambda'\brac{r} \brac{g_{\Sigma}}_{ij} - u\brac{\II_{\Sigma}}_{ij} }.\label{fml:HessLambda'}
\end{equation}

The following Minkowski type formula can be found in \cite{LWX26}*{Lem.~2.5}; see also \cite{HLW22}*{Lem.~2.6}.
\begin{lem}
Let $\varepsilon \in \R$, and let $\Sigma$ be a closed $C^2$ hypersurface in $\M^{n+1}$. Then
\begin{equation} \label{eq:shiftedMinForm}
\int_{\Sigma} \brac{\lambda'\brac{r}- \varepsilon u} p_{m-1} \brac{\kappa -\varepsilon} \, d\mu = \int_{\Sigma} u p_{m} \brac{\kappa - \varepsilon} \, d\mu, \quad m=1, \ldots, n.
\end{equation}
\end{lem}
For $k=0,1,\ldots, n$, equation \eqref{eq:diverW0} and the above Minkowski type formula \eqref{eq:shiftedMinForm} yield the following expression for the $k$th weighted shifted curvature integral defined in \eqref{eq:WeiCurInte}:
\begin{align}\label{Wk-EquiDef}
W_k^{\lambda', \varepsilon} \brac{\Omega} = \int_{\Sigma} u p_k \brac{\kappa -\varepsilon} \, d\mu.
\end{align}

\begin{prop}\label{prop:Wkeps-Equidef-Wk}
Let $\varepsilon \in \R$, and let $\Omega$ be a bounded domain in $\M^{n+1}$ with $C^2$ boundary. Then
\begin{equation}\label{eq:Wkeps-Wk}
W_k^{\lambda',\varepsilon} \brac{\Omega} = \sum_{i=0}^k \binom{k}{i}  \brac{-\varepsilon}^{k-i} W_i^{\lambda'} \brac{\Omega}, \quad k=0,\ldots, n+1.
\end{equation}
\end{prop}
\begin{proof}
First, we assume that $0\leq k \leq n$. Then by \eqref{Wk-EquiDef},
\begin{align}
W_k^{\lambda',\varepsilon} \brac{\Omega} &= \int_{\Sigma} u p_k \brac{\kappa - \varepsilon} \, d\mu = \int_{\Sigma} u \brac{\sum_{i=0}^k \binom{k}{i} p_i \brac{\kappa} \brac{-\varepsilon}^{k-i} } \, d\mu \nonumber\\
&\qquad = \sum_{i=0}^k \binom{k}{i} \brac{-\varepsilon}^{k-i} W_i^{\lambda'} \brac{\Omega}. \label{eq:WkL-0-n}
\end{align}
Using the case $k=n$ of \eqref{eq:WkL-0-n}, we have
\begin{align*}
W_{n+1}^{\lambda',\varepsilon} \brac{\Omega} &= \int_{\Sigma} \brac{\lambda'\brac{r} - \varepsilon u} p_n \brac{\kappa -\varepsilon} \, d\mu\\
&= \int_{\Sigma} \lambda' \brac{r} \brac{\sum_{i=0}^n \binom{n}{i} p_i \brac{\kappa} \brac{-\varepsilon}^{n-i} } \, d\mu
-\varepsilon W_n^{\lambda',\varepsilon} \brac{\Omega}\\
&=\sum_{j=1}^{n+1} \binom{n}{j-1} \brac{-\varepsilon}^{n-j+1} W_j^{\lambda'} \brac{\Omega}+ \sum_{i=0}^n \binom{n}{i} \brac{-\varepsilon}^{n-i+1} W_i^{\lambda'} \brac{\Omega}\\
&= \sum_{j=0}^{n+1} \binom{n+1}{j} \brac{-\varepsilon}^{n-j+1} W_j^{\lambda'} \brac{\Omega}.
\end{align*}
This completes the proof.
\end{proof}

The following corollary of Proposition \ref{prop:Wkeps-Equidef-Wk} was observed in \cite{LWX26}*{Sect.~4}. 
\begin{cor}\label{cor-LWXObser}
Let $\varepsilon \in \R$, and let $\Omega$ be a bounded domain in $\M^{n+1}$ with $C^2$ boundary. Then
\begin{equation*}
W_1^{\lambda', \varepsilon} \brac{\Omega}^2 - W_0^{\lambda',\varepsilon} \brac{\Omega} W_2^{\lambda', \varepsilon} \brac{\Omega} = W_1^{\lambda'}\brac{\Omega}^2 -W_0^{\lambda'} \brac{\Omega} W_2^{\lambda'}\brac{\Omega}.
\end{equation*}
\end{cor}
\begin{proof}
Using \eqref{eq:Wkeps-Wk}, we have
\begin{align*}
W_0^{\lambda',\varepsilon} \brac{\Omega} &= W_0^{\lambda'} \brac{\Omega}, \quad W_1^{\lambda',\varepsilon} \brac{\Omega} = W_1^{\lambda'} \brac{\Omega} - \varepsilon W_0^{\lambda'} \brac{\Omega},\\
W_2^{\lambda',\varepsilon} \brac{\Omega} &=W_2^{\lambda'} \brac{\Omega} - 2\varepsilon W_1^{\lambda'}\brac{\Omega}+ \varepsilon^2 W_0^{\lambda'} \brac{\Omega}.
\end{align*}
The assertion then follows from direct calculation.
\end{proof}

A domain $\Omega$ (or its boundary $\Sigma$) in $\M^{n+1}$ is called static convex if 
\begin{equation*}
\II_{\Sigma} \geq \frac{\lambda'' \brac{r} }{\lambda'\brac{r}}\sscalar{\partial_r, \nu} g_{\Sigma}
\end{equation*}
on $\Sigma$, and it is called strictly static convex if the inequality above is strict. If $\M^{n+1} = \M^{n+1} \brac{c}$, we define 
\begin{equation}\label{def:S-Sigma}
S_{\Sigma} = \II_{\Sigma} + c \frac{u_{\Sigma}}{\lambda'\brac{r}} g_{\Sigma}.
\end{equation}
Since 
\begin{equation*}
\lambda''\brac{r}\sscalar{ \partial_r, \nu}/ \lambda'\brac{r} = -c \sscalar{\lambda\brac{r} \partial_r, \nu} / \lambda'\brac{r} = -c u_{\Sigma}/\lambda'\brac{r},
\end{equation*}
$\Sigma$ is static convex if and only if $S_{\Sigma}$ is nonnegative definite on $\Sigma$. By \eqref{def:S-Sigma}, static convexity reduces to the usual notion of convexity when $\M^{n+1} = \R^{n+1}$.

In the following, we verify that all geodesic balls in $\M^{n+1}$ are strictly static convex. This fact has been implicitly used in the literature, although an explicit proof seems to be absent.
\begin{prop}\label{prop:BallStaticConvex}
Let $\Omega$ be a geodesic ball in $\M^{n+1} \in \set{\H^{n+1}, \S^{n+1}_+}$. Then $\Omega$ is strictly static convex.
\end{prop}
\begin{proof}
\textbf{Case 1. $\M^{n+1} = \H^{n+1}$.} Let $\Omega \subset \H^{n+1}$ be a geodesic ball of radius $\rho$. Then at any $X \in \Sigma$, we have $\abs{u_{\Sigma}/ \lambda'\brac{r}} \leq \lambda\brac{r}/ \lambda'\brac{r} <1$. Meanwhile, we have $\II_{\Sigma} = \lambda'\brac{\rho} / \lambda \brac{\rho} g_{\Sigma}>g_{\Sigma}$. Thus, $S_{\Sigma}$ is strictly positive definite on $\Sigma$. 

\textbf{Case 2. $\M^{n+1} = \S^{n+1}_+$.} Let $\Omega \subset \S^{n+1}_+$ be a geodesic ball of radius $\rho$ centered at $P_0$. By \eqref{eq:geoS}, for any $P \in \Sigma$, there exists $v \in T_{P_0} \S^{n+1}$ such that $P = \lambda'\brac{\rho} P_0 + \lambda\brac{\rho} v$. Consider the geodesic $\gamma\brac{t} = \lambda'\brac{t} P_0 + \lambda\brac{t} v $ from $P_0$ to $P$, where $t \in [0, \rho]$. Then we have at $P$,
\begin{equation*}
\partial_{\nu} \brac{\log \lambda' \brac{r}} =\left. \frac{d}{dt} \right|_{t=\rho} \log \sscalar{\gamma \brac{t}, (0,1)} = \frac{- \lambda \brac{\rho}\sscalar{P_0, (0,1)} + \lambda'\brac{\rho} \sscalar{v, (0,1)}   }{\lambda'\brac{\rho} \sscalar{P_0, (0,1)} + \lambda\brac{\rho} \sscalar{v, (0,1)} }.
\end{equation*}
Since the principal curvatures of $\Sigma$ at $P$ are all equal to $\lambda'\brac{\rho} / \lambda \brac{\rho}$, we can use the above formula to obtain 
\begin{align*}
S_{\Sigma} = \II_{\Sigma} - \partial_{\nu}\brac{\log \lambda' \brac{r} }g_{\Sigma} = \frac{\sscalar{P_0, (0,1)}}{\lambda\brac{\rho} \lambda'\brac{r(P)}} g_{\Sigma}>0.
\end{align*}
This shows that $\Omega$ is strictly static convex.
\end{proof}

\subsection{Approximation lemmas}\label{subsec:approx}
Let $\Sigma$ be a closed $C^3$ static convex hypersurface in $\M^{n+1}(c)$ satisfying $\II_{\Sigma} \geq \varepsilon g_{\Sigma}$. In this subsection, we prove that $\Sigma$ can be approximated by a sequence of strictly static convex hypersurfaces with principal curvatures strictly greater than $\varepsilon$. 

The following evolution equations can be found in \cite{LWX26}*{Lem.~2.6}.
\begin{lem}
Let $X: \Sigma \times [0, T) \to \M^{n+1} (c)$ be a family of immersions with $\Sigma_t= X \brac{\Sigma, t}$ evolving by $\partial_t X = \mathcal{F} \nu$. Then the following evolution equations hold:
\begin{align}
\partial_t \lambda' \brac{r} &=  -cu \mathcal{F}, \label{fml:evl-Lam}  \\ 
\partial_t u &= \lambda' \brac{r} \mathcal{F}-\sscalar{\lambda \brac{r} \partial_r, \nabla \mathcal{F}}, \label{fml:evl-U}\\
\partial_t \II_i{}^j &= - \nabla^j \nabla_i \mathcal{F} - \mathcal{F} \brac{\brac{\II^2}_i{}^j +c \delta_i{}^j}, \label{fml:evl-II}
\end{align}
where $\II$ is short for $\II_{\Sigma_t}$, and $u$ is short for $u_{\Sigma_t}$.
\end{lem}

\begin{lem}
Let $C \in \R$ be a constant. Along the flow $\partial_t X =\mathcal{F} \nu$ in $\M^{n+1} \brac{c}$ with $\mathcal{F} = \lambda'\brac{r} -C$, the following evolution equations hold:
\begin{align}
\partial_t \II_i{}^j &=  C\brac{\brac{\II^2}_i{}^j + c \delta_i{}^j}-\lambda'\brac{r} \II_i{}^k S_k{}^j, \label{evl-2FF-rea}\\
\partial_t S_i{}^j &=\brac{ \brac{C-\lambda'\brac{r}} \II_i{}^k- c C\lambda'\brac{r}^{-1} u \delta_i{}^k } S_k{}^j + c \lambda'\brac{r}^{-1} \delta_i{}^j, \label{evl-S-rea}
\end{align}
where $S_i{}^j$ is short for $\brac{S_{\Sigma_t} g_{\Sigma_t}^{-1}}_i{}^j$.
\end{lem}

\begin{proof}
Using \eqref{fml:evl-II} and \eqref{fml:HessLambda'}, we have
\begin{align*}
\partial_t \II_i{}^j &= - \nabla^j \nabla_i \brac{\lambda' \brac{r} -C} - \brac{\lambda'\brac{r} -C} \brac{ \brac{\II^2}_i{}^j + c \delta_i{}^j} \\
&= c \brac{\lambda'\brac{r} \delta_i{}^j - u \II_i{}^j}- \brac{\lambda'\brac{r} -C} \brac{ \brac{\II^2}_i{}^j + c \delta_i{}^j} \\
&=C \brac{\brac{\II^2}_i{}^j + c \delta_i{}^j} - \II_i{}^{k} \brac{\lambda'\brac{r} \II_k{}^j + c u \delta_k{}^j}\\
&= C\brac{\brac{\II^2}_i{}^j + c \delta_i{}^j}-\lambda'\brac{r} \II_i{}^k S_k{}^j.
\end{align*}
Thus, we obtain \eqref{evl-2FF-rea}.

Since $\lambda'' = - c \lambda$, we know that $\nabla \lambda' \brac{r}$ is the projection of the vector field $-c \lambda \brac{r} \partial_r$ onto $T_X \Sigma$. Thus we have $\sscalar{\lambda \brac{r} \partial_r, \nabla \lambda' \brac{r}} = -c \lambda\brac{r}^2 \brac{1- \sscalar{\partial_r, \nu}^2}= -1+\lambda'\brac{r}^2 + cu^2$, where we used $(\lambda')^2+c\lambda^2=1$ in the last equality. Then by \eqref{fml:evl-U}, we have
\begin{align*}
\partial_t u &= \lambda'\brac{r} \brac{\lambda'\brac{r} -C} - \sscalar{\lambda\brac{r} \partial_r, \nabla \brac{\lambda'\brac{r} -C }}\\
&= \lambda'\brac{r}^2- C \lambda'\brac{r}+1-\lambda'\brac{r}^2-cu^2\\
&=1-C \lambda'\brac{r} -cu^2.
\end{align*}
Using the above formula and \eqref{fml:evl-Lam}, we calculate 
\begin{align*}
\partial_t \frac{u}{\lambda'\brac{r}} &= \lambda'\brac{r}^{-1} -C - c \lambda'\brac{r}^{-1} u^2+ c\lambda'\brac{r}^{-2} u^2 \brac{\lambda'\brac{r} -C}\\
&= \lambda'\brac{r}^{-1} - C \brac{1+c \lambda'\brac{r}^{-2} u^2}.
\end{align*}
Combining the above evolution equation with \eqref{evl-2FF-rea} yields
\begin{align*}
\partial_t S_i{}^j&= C\brac{\brac{\II^2}_i{}^j + c \delta_i{}^j}-\lambda'\brac{r} \II_i{}^k S_k{}^j+c \brac{\lambda'\brac{r}^{-1} - C \brac{1+c \lambda'\brac{r}^{-2} u^2} } \delta_i{}^j\\
&=C \brac{\brac{\II^2}_i{}^j - c^2\lambda'\brac{r}^{-2}u^2 \delta_i{}^j }- \lambda'\brac{r} \II_i{}^k S_k{}^j + c \lambda'\brac{r}^{-1} \delta_i{}^j\\
&= \brac{ \brac{C-\lambda'\brac{r}} \II_i{}^k- c C\lambda'\brac{r}^{-1} u \delta_i{}^k } S_k{}^j + c \lambda'\brac{r}^{-1} \delta_i{}^j.
\end{align*}
Hence, we obtain \eqref{evl-S-rea}. This completes the proof.
\end{proof}

In the following Lemma \ref{lem:approximation-H} and Lemma \ref{lem:approximation-S}, we establish the desired approximation by constructing a flow starting from $\Sigma$ with a suitable normal velocity.
\begin{lem}\label{lem:approximation-H}
Let $\Omega$ be a bounded domain in $\H^{n+1}$ with $C^3$ boundary $\Sigma$. Let $\varepsilon \in \R$ be a constant. Assume that $\Sigma$ is static convex with $\II_{\Sigma} \geq \varepsilon g_{\Sigma}$. Then there exists a sequence of closed $C^2$ hypersurfaces $\Sigma_i$ such that each $\Sigma_i$ is strictly static convex with $\II_{\Sigma_i} > \varepsilon g_{\Sigma_i}$, and $\Sigma_i \to \Sigma$ in $C^2$-topology as $i \to \infty$. 
\end{lem}

\begin{proof}
For any $f \in C^2\brac{\Sigma}$, there exists $\delta>0$ such that for any $t \in [0,\delta)$,
\begin{equation*}
\Sigma_t = \set{ \lambda' \brac{ f \brac{X} t} X+  \lambda \brac{f \brac{X} t}  \nu \brac{X};~X \in \Sigma}
\end{equation*}
is a closed $C^2$ hypersurface in $\H^{n+1}$. Since $\lambda'\brac{0} =1$, the hypersurfaces $\Sigma_t$ define a deformation of $\Sigma$ with initial velocity $\mathcal{F} \nu$, where $\mathcal{F} = f(X)$. In the following, we construct specific $\mathcal{F}$ such that $\Sigma_t$ is strictly static convex with $\II_{\Sigma_t} > \varepsilon g_{\Sigma_t}$ for small $t>0$, which yields the desired approximating hypersurfaces $\Sigma_i$.  For simplicity, we write $\II_{\Sigma_t}$ as $\II$, $u_{\Sigma_t}\brac{X_t}$ as $u$, and $\lambda' \brac{r \brac{X_t}}$ as $\lambda'\brac{r}$.

\textbf{Case 1. $\varepsilon \in (-\infty,1)$.} In this case, we take $\mathcal{F} = C - \lambda'\brac{r}$, where $C$ is a positive constant. By \eqref{evl-2FF-rea}, we have
\begin{align}
\partial_t \II_i{}^j = C \brac{\delta_i{}^j - \brac{\II^2}_i{}^j } + \lambda'\brac{r} \II_i{}^k S_k{}^j. \label{fml:VarIIReal}
\end{align}
If $\varepsilon \in (-1,1)$, we can choose $C$ sufficiently large so that $\partial_t \II_i{}^j$ is strictly positive in every null direction of $\II_i{}^j-\varepsilon \delta_i{}^j$ at $t=0$. This implies $\II_{\Sigma_t} > \varepsilon g_{\Sigma_t}$ for small $t>0$. If $\varepsilon \in (-\infty, -1]$, the static convexity of $\Sigma$ directly implies
\begin{equation*}
\II_{\Sigma} \geq \frac{u}{\lambda'\brac{r}} g_{\Sigma} \geq -\frac{\lambda\brac{r}}{\lambda'\brac{r}} g_{\Sigma} > -g_{\Sigma}.
\end{equation*}
Thus $\II_{\Sigma} > \varepsilon g_{\Sigma}$, and this property is stable under perturbation.

On the other hand, by \eqref{evl-S-rea} we have
\begin{equation*}
\partial_t S_i{}^j = -\brac{\brac{C-\lambda'} \II_i{}^k + C  \lambda'\brac{r}^{-1} u \delta_i{}^k} S_k{}^j+\brac{\lambda'\brac{r}}^{-1} \delta_i{}^j.
\end{equation*}
Since the right-hand side is strictly positive in any null direction of $\brac{S_i{}^j}$, we have that $\Sigma_t$ is strictly static convex for small $t>0$. 

\textbf{Case 2. $\varepsilon \in [1, \infty)$.} In this case, we take $\mathcal{F} = -\lambda'\brac{r}$. Since $\II_{\Sigma} \geq \varepsilon g_{\Sigma} > \abs{u/\lambda'(r)}g_{\Sigma}$, we have that $S_i{}^j$ is strictly positive on $\Sigma_t$ for small $t \geq 0$. Then by taking $C=0$ in \eqref{fml:VarIIReal}, we have at $t=0$, 
\begin{equation*}
\partial_t \II_i{}^j = \lambda'\brac{r} \II_i{}^k S_k{}^j>0.
\end{equation*}
Thus $\II_{\Sigma_t} > \varepsilon g_{\Sigma_t}$ for small $t>0$.
This completes the proof.
\end{proof}

\begin{lem}\label{lem:approximation-S}
Let $\Omega$ be a bounded domain in $\S^{n+1}_+$ with $C^3$ boundary $\Sigma$. Let $\varepsilon \in \R$ be a constant. Assume that $\Sigma$ is static convex with $\II_{\Sigma} \geq \varepsilon g_{\Sigma}$. Then there exists a sequence of closed $C^2$ hypersurfaces $\Sigma_i$ such that each $\Sigma_i$ is strictly static convex with $\II_{\Sigma_i}>\varepsilon g_{\Sigma_i}$, and $\Sigma_i \to \Sigma$ in $C^2$-topology as $i \to \infty$. 
\end{lem}

\begin{proof}
We consider the flow $\partial_t X = \mathcal{F} \nu$ with initial speed $\mathcal{F} = \lambda' \brac{r}-C$, where $C$ is a large positive constant so that $C> \sup_{X \in \Sigma} \norm{\II \circ S}_{op}$. Then we have by \eqref{evl-2FF-rea} that $\II_{\Sigma_t} >  \varepsilon g_{\Sigma_t}$ for small $t>0$. On the other hand, it follows by \eqref{evl-S-rea} that $\Sigma_t$ is strictly static convex for small $t>0$. The assertion follows.
\end{proof}

\section{Orthogonal projection} \label{sec:Gib}
To study a conjectured inequality of Penrose for spacelike $2$-surfaces in $\R^{3,1}$, Gibbons \cite{Gib97} proposed an approach based on orthogonal projection onto a Euclidean hyperplane and the Minkowski inequality for the projected surface. This approach is known as Gibbons' projection procedure.
Later, Brendle and Wang \cite{BW14} extended this approach to Schwarzschild spacetime and proved the corresponding Gibbons–Penrose inequality.

Now we project the hyperboloid $\H^{n+1}$ orthogonally onto $\R^{n+1} \times \set{0}$ in $\R^{n+1,1}$ and define  $\pi_{\H} : \H^{n+1} \to \R^{n+1}$ by $\pi_{\H} \brac{x, x_{n+2}} = x$. This projection is a diffeomorphism. Similarly, we define $\pi_{\S} : \S^{n+1}_+ \to B^{n+1}$ by $\pi_{\S} \brac{p, p_{n+2}} =p$, where $B^{n+1}$ denotes the centered unit ball in $\R^{n+1} \times \set{0}$; see Figure \ref{Fig:GibProj}. When necessary, we also extend $\pi_{\H}$ to $\R^{n+1,1}$ by $\pi_{\H} \brac{x, x_{n+2}}=x$ and $\pi_{\S}$ to $\R^{n+2}$ by $\pi_{\S} \brac{p, p_{n+2}} =p$. For simplicity, we denote both $\pi_{\H}$ and $\pi_{\S}$ by $\pi$ whenever no confusion arises. 
\begin{figure}[htbp]
\centering
\begin{tikzpicture} 
\coordinate (X) at (1, {sqrt(2)});
\draw[->] (-1.7,0) -- (2.5,0) node[right]{$x$};
\draw[->] (0,-0.2)-- (0,2.5) node[above]{$x_{n+2}$};
\draw[domain= -1:2, samples=200, thick]
plot(\x, {sqrt(\x*\x+1)});
\foreach \x in {-0.7,-0.133333,0.433333,1} {
\draw[thick, dashed, ->]
(\x, {sqrt(\x*\x+1)}) -- (\x, 0);
}
\node at (-1.4,1.5){$\mathbb{H}^{n+1}$};
\node[above=4pt] at (X) {$X$};
\node at (1.4, -0.4){$\pi_{\H} \brac{X}$};
\node at (1.4,0.7){$\pi_{\H}$};
\filldraw (X) circle(1.2pt);
\end{tikzpicture}
\qquad
\qquad
\begin{tikzpicture} 
\coordinate (P) at (0.6, 1.42);
\draw[->] (-1.7,0) -- (2.5,0) node[right]{$p$};
\draw[->] (0,-0.2)-- (0,2.5) node[above]{$p_{n+2}$};
\draw[thick] (1.5,0.05)
  arc[start angle=0, end angle=180, radius=1.5];
\foreach \x in {-0.8,-0.34,0.14,0.6} {
\draw[thick, dashed, ->]
(\x, {sqrt(2.25-\x*\x)+0.05}) -- (\x, 0);
}
\node at (-1.4,1.5){$\mathbb{S}^{n+1}_+$};
\node[above=4pt] at (P) {$P$};
\node[below=4pt] at (1,0) {$\pi_{\S} \brac{P}$};
\node at (0.9,0.7){$\pi_{\S}$};
\filldraw (P) circle(1.2pt);
\end{tikzpicture}
\caption{Orthogonal projections}
\label{Fig:GibProj}
\end{figure}

In this section, we assume that either $\M^{n+1} = \H^{n+1}$ or $\M^{n+1}= \S^{n+1}_+$. The following lemma was proved in \cite{LX23}.
\begin{lem}[\cite{LX23}]
Let $\Omega$ be a bounded domain in $\M^{n+1}$. Then
\begin{equation}\label{eq:Vol-WeightedVol}
W_0^{\lambda'} \brac{\Omega} = \brac{n+1} \Vol \brac{\pi \brac{\Omega}}.
\end{equation}
\end{lem}

\begin{proof}
Since the projection $\pi$ maps the radial coordinate $r$ of $\M^{n+1}$ to the radial coordinate $\rho = \lambda \brac{r}$ of $\R^{n+1}$, we have $d\rho = \lambda'\brac{r} dr$. Let $d\vol_{\M}$ and $d\vol_{\R}$ denote the volume elements of $\M^{n+1}$ and $\R^{n+1}$ respectively, and let $d\sigma$ denote the area element of the unit sphere $\S^n$. Then
\begin{equation}\label{eq:rel-VolumeMeasure}
\pi^* \brac{d\vol_{\R}} = \pi^* \brac{\rho^n d\rho \, d\sigma} = \lambda\brac{r}^n \lambda' \brac{r}\, dr \, d\sigma = \lambda' \brac{r} \, d\vol_{\M}.
\end{equation}
The lemma follows.
\end{proof}

Throughout this paper, for any hypersurface $\Sigma$ in $\M^{n+1}$, we let $\widehat{\Sigma} = \pi \brac{\Sigma}$ be its projection in $\R^{n+1}$. Using the same local coordinates for $\Sigma$ and $\widehat{\Sigma}$ via the map $\pi$, we can compare their geometries. For the case $\M^{n+1} = \H^{n+1}$, this was studied in \cites{HLWX25, LX23}. We begin by recalling some identities from these papers. The case $\M^{n+1} =\S^{n+1}_+$ is addressed separately in Subsection \ref{subsec:GibS}.

\subsection{Hyperbolic case}
\begin{lem}[\cite{LX23}]\label{lem:nu(X)}
Let $\Omega$ be a bounded domain in $\H^{n+1}$ with $C^1$ boundary $\Sigma$. Then at any $X= \brac{x, x_{n+2}} \in \Sigma$, we have
\begin{equation}\label{eq:GibH-nu}
\nu_{\Sigma} \brac{X} = \frac{\brac{\nu_{\widehat{\Sigma}} + u_{\widehat{\Sigma}} x, x_{n+2}u_{\widehat{\Sigma}}   }}{\sqrt{1+u_{\widehat{\Sigma}}^2 } },
\end{equation} 
where $\nu_{\widehat{\Sigma}}$ is the outward unit normal to $\widehat{\Sigma}$ at $x$, and $u_{\widehat{\Sigma}}$ is the support function of $\widehat{\Sigma}$ at $x$.
\end{lem}

\begin{lem}[\cite{HLWX25}]\label{lem:GibH-GeoQuantitu}
Let $\Sigma$ be a $C^2$ hypersurface in $\H^{n+1}$. Then
\begin{align}
\lambda' \brac{r} &=\sqrt{1+ \abs{x}^2}, \label{eq:GibH-lambda'}\\
\brac{g_{\Sigma}}_{ij} &= \brac{g_{\widehat{\Sigma}}}_{ij} - \frac{ \sscalar{x, \partial_i x} \sscalar{x, \partial_j x} }{1+ \abs{x}^2}, \label{eq:GibH-metric}\\
u_{\Sigma} &= \brac{\frac{1+ \abs{x}^2}{1+ u_{\widehat{\Sigma}}^2} }^{\frac{1}{2}} u_{\widehat{\Sigma}}, \label{eq:GibH-support}\\
\brac{\II_{\Sigma}}_{ij} &= \frac{1}{\sqrt{1+u_{\widehat{\Sigma}}^2}} \brac{ \brac{\II_{\widehat{\Sigma}}}_{ij} +u_{\widehat{\Sigma}} \brac{g_{\widehat{\Sigma}}}_{ij} - \frac{u_{\widehat{\Sigma}} \sscalar{x, \partial_i x} \sscalar{x, \partial_j x} }{1+ \abs{x}^2} }, \label{fml-2FFSigma}
\end{align}
where the left-hand sides are evaluated at $X = \brac{x, x_{n+2}} \in \Sigma$ and the right-hand sides are evaluated at $x \in \widehat{\Sigma}$.
\end{lem}

Taking $c=-1$ in \eqref{def:S-Sigma} and using Lemma \ref{lem:GibH-GeoQuantitu}, we have
\begin{equation*}
\lambda'\brac{r} \brac{S_{\Sigma}}_{ij} =\lambda'\brac{r}  \brac{\II_{\Sigma}}_{ij}-u _{\Sigma} \brac{g_{\Sigma}}_{ij} = \brac{\frac{1+\abs{x}^2}{1+ u_{\widehat{\Sigma}}^2} }^{\frac{1}{2}} \brac{\II_{\widehat{\Sigma}}}_{ij}.
\end{equation*}
This implies the following proposition; see \cite{HLWX25}*{p.~69}.
\begin{prop}\label{prop:GibH-StaticC-Convex}
Let $\Sigma$ be a $C^2$ hypersurface in $\H^{n+1}$. Then $\Sigma$ is (strictly) static convex in $\H^{n+1}$ if and only if $\widehat{\Sigma}$ is (strictly) convex in $\R^{n+1}$.
\end{prop}

\begin{proof}[Proof of statements \eqref{prop3.1} and \eqref{prop3.2} of Proposition \ref{prop:ConvComp}]
Let $\Omega \subset \H^{n+1}$ be a bounded static convex domain with $C^2$ boundary $\Sigma$. By the static convexity condition \eqref{cond:StaticConvex} and $\abs{\lambda''\brac{r}/\lambda'\brac{r}} = \abs{\lambda\brac{r}/\lambda'\brac{r}}<1$ on $\Sigma$, we have $\II_{\Sigma} >- g_{\Sigma}$. This proves statement \eqref{prop3.1}. 

If in addition $\mathcal{O} \in \mathrm{int}\brac{\Omega}$, then by Proposition \ref{prop:GibH-StaticC-Convex}, 
$\pi \brac{\Omega} \subset \R^{n+1}$ is convex and contains the origin in its interior. It follows that $u_{\widehat{\Sigma}}>0$, and thus $u_{\Sigma}>0$ by \eqref{eq:GibH-support}. The static convexity condition \eqref{cond:StaticConvex} then yields $\II_{\Sigma}>0$. This proves statement \eqref{prop3.2}.
\end{proof}

Let $\Omega \subset \H^{n+1}$ be a strictly static convex domain with boundary $\Sigma$. By Proposition \ref{prop:GibH-StaticC-Convex}, $K_{\Omega} := \pi \brac{\Omega}$ is strictly convex in $\R^{n+1}$. We can therefore use its support function $h_{K_{\Omega}}$ to characterize the geometry of $\partial K_{\Omega}$. In the following lemma, we translate the condition $\II_{\Sigma} > \varepsilon g_{\Sigma}$ on $\Sigma$ into a condition on $h_{K_{\Omega}}$. 

\begin{lem}\label{lem:Convex-Tensor-H}
Let $\Omega$ be a bounded domain in $\H^{n+1}$ with $C^2$ boundary $\Sigma$, and let $K_{\Omega} = \pi \brac{\Omega}$. Assume that $\Sigma$ is strictly static convex and that $\II_{\Sigma} > \varepsilon g_{\Sigma}$ for some $\varepsilon \in \R$. Then the tensor
\begin{equation}\label{TenSConv-H}
\brac{{h_{K_{\Omega}}}\brac{1+h_{K_{\Omega}}^2}^{-\frac{1}{2}}-\varepsilon}D^2 h_{K_{\Omega}} + \brac{ 1+h_{K_{\Omega}}^2}^{-\frac{3}{2}} D h_{K_{\Omega}} \otimes D h_{K_{\Omega}}  + \brac{\brac{1+h_{K_{\Omega}}^2}^{\frac{1}{2}} - \varepsilon h_{K_{\Omega} } } g_{\S^n}
\end{equation}
is positive definite on $\S^n$.
\end{lem}

\begin{proof}
Using \eqref{fml-2FFSigma} and \eqref{eq:GibH-metric}, we have
\begin{equation}\label{eq:II-metric-GibH}
\brac{\II_{\Sigma}- \varepsilon g_{\Sigma}}_{ij} =  \frac{1}{\sqrt{1+u_{\widehat{\Sigma}}^2}} \brac{\II_{\widehat{\Sigma}}}_{ij} + \brac{\frac{u_{\widehat{\Sigma}} }{ \sqrt{1+u_{\widehat{\Sigma}}^2} }- \varepsilon } \brac{\brac{g_{\widehat{\Sigma}}}_{ij}- \frac{\sscalar{x, \partial_i x} \sscalar{x, \partial_j x} }{1+ \abs{x}^2}}.
\end{equation}
By the assumption $\II_{\Sigma}> \varepsilon g_{\Sigma}$, the right-hand side of \eqref{eq:II-metric-GibH} is positive definite on $\widehat{\Sigma} = \partial K_{\Omega}$.

Write $h$ for $h_{K_{\Omega}}$ as an abbreviation. Since $\Omega$ is strictly static convex, Proposition \ref{prop:GibH-StaticC-Convex} implies that $K_{\Omega}$ is strictly convex, and thus $D^2 h+ h g_{\S^n}$ is positive definite on $\S^n$. If $h\brac{1+h^2}^{-\frac{1}{2}} - \varepsilon \geq 0$ at $\theta \in \S^n$, then
\begin{align*}
&\brac{h \brac{1+h^2}^{-\frac{1}{2}}-\varepsilon} D^2 h + \brac{1+h^2}^{-\frac{3}{2}} Dh \otimes Dh+ \brac{\brac{1+h^2}^{\frac{1}{2}}-\varepsilon h} g_{\S^n} \\
&\geq -h \brac{h \brac{1+h^2}^{-\frac{1}{2}}-\varepsilon} g_{\S^n} + \brac{\brac{1+h^2}^{\frac{1}{2}}-\varepsilon h} g_{\S^n}= \brac{1+h^2}^{-\frac{1}{2}} g_{\S^n}>0.
\end{align*}
Hence, it remains to consider the points where $h- \varepsilon \brac{1+h^2}^{\frac{1}{2}} <0$.

Assume that $h- \varepsilon \brac{1+h^2}^{\frac{1}{2}} <0$ at $\theta \in \S^n$. Let $x = D h \brac{\theta} +h \theta$ be the unique point on $\widehat{\Sigma}$ with the outward unit normal $\theta$. Take an orthonormal basis $\set{e_1, \ldots, e_n}$ in $T_x \widehat{\Sigma}$. Note that it is also an orthonormal basis for $T_\theta \S^n$. By \eqref{eq:II-metric-GibH}, \eqref{eq:Rel-Supports}, and \eqref{fml:2FF-D2h}, the matrix 
\begin{equation*}
\brac{D^2 h + h \mathrm{Id}}^{-1}+ \brac{h-\varepsilon \brac{1+h^2}^{\frac{1}{2}} }\brac{ \mathrm{Id}-\frac{Dh \otimes Dh }{1+ \abs{D h}^2+ h^2}}
\end{equation*}
is positive definite at $\theta \in \S^n$. Equivalently,
\begin{equation*}
\brac{D^2 h+ h \mathrm{Id}}^{-1} > \brac{\varepsilon \brac{1+h^2}^{\frac{1}{2}}-h } \brac{\mathrm{Id} -\frac{Dh \otimes Dh }{1+ \abs{D h}^2+ h^2} }.
\end{equation*}
Taking the inverse of both sides yields 
\begin{equation*}
D^2 h+ h \mathrm{Id} < \frac{1}{\varepsilon \brac{1+h^2}^{\frac{1}{2}}-h  } \brac{\mathrm{Id} +\frac{Dh \otimes Dh}{1+ h^2}}.
\end{equation*}
Rearranging this inequality gives
\begin{equation*}
\brac{\varepsilon \brac{1+h^2}^{\frac{1}{2} }-h} D^2h + \brac{\varepsilon h\brac{1+h^2}^{\frac{1}{2}}-h^2-1}\mathrm{Id}
-\frac{Dh \otimes Dh}{1+h^2}<0.
\end{equation*}
Multiplying both sides by $-(1+h^2)^{-\frac12}$, we obtain
\begin{equation*}
\brac{h\brac{1+h^2}^{-\frac{1}{2}}-\varepsilon}D^2h +\brac{1+h^2}^{-\frac{3}{2}}Dh \otimes Dh +\brac{\brac{1+h^2}^{\frac{1}{2}} -\varepsilon h} \mathrm{Id}>0.
\end{equation*}
This completes the proof.
\end{proof}

In the following lemma, the definition of $\pi$ is extended to $\R^{n+1,1}$ by $\pi \brac{x, x_{n+2}}= x$.
\begin{lem}\label{lem:regularity-Lomega-H}
Let $\Omega$ be a bounded domain in $\H^{n+1}$ with $C^3$ boundary $\Sigma$. Assume that $\II_{\Sigma} > \varepsilon g_{\Sigma}$ for some $\varepsilon \in \R$. Then the map $F : \Sigma \to \R^{n+1}$ defined by $F (X) = \pi \brac{\nu \brac{X} -\varepsilon X}$ is a $C^2$ immersion.
\end{lem}

\begin{proof}
Define $\bar{F}: \Sigma \to \R^{n+1,1}$ by $\bar{F} \brac{X} = \nu \brac{X} -\varepsilon X$. Since $\II_{\Sigma} > \varepsilon g_{\Sigma}$, $d\bar{F}$ has rank $n$. Since for any $0 \neq v \in T_X \Sigma$, $d\bar{F} \brac{v} \in T_X \Sigma$ is a nonzero spacelike vector, and $\ker \brac{d\pi} = \R \brac{0,1}$ is timelike, we have that $d \brac{\pi \circ \bar{F}}$ has rank $n$. Hence, $F = \pi \circ \bar{F}$ is an immersion of $\Sigma$ into $\R^{n+1}$. 
\end{proof}

\subsection{Spherical case}\label{subsec:GibS}
\begin{lem}\label{lem:GibS-GeoQuantitu}
Let $\Sigma$ be a $C^2$ hypersurface in $\S^{n+1}_+$. Then
\begin{align}
\lambda' \brac{r} &= p_{n+2} =\sqrt{1- \abs{p}^2},\label{GibS-lambda'}\\
\nu_{\Sigma} &= \frac{\brac{\nu_{\widehat{\Sigma}} - u_{\widehat{\Sigma}} p, - p_{n+2} u_{\widehat{\Sigma}}} }{\sqrt{1-u_{\widehat{\Sigma}}^2}} \label{GibS-nu},\\
\brac{g_{\Sigma}}_{ij} &= \brac{g_{\widehat{\Sigma}}}_{ij}+ \frac{\sscalar{p,\partial_i p}\sscalar{p, \partial_j p} }{1-\abs{p}^2},\label{GibS-metric}\\
u_{\Sigma} &= \brac{\frac{1-\abs{p}^2 }{1-u_{\widehat{\Sigma}}^2}}^{\frac{1}{2}} u_{\widehat{\Sigma}}, \label{GibS-supp}\\
\brac{\II_{\Sigma}}_{ij} &=\frac{1}{\sqrt{1-u_{\widehat{\Sigma}}^2}} \brac{\brac{\II_{\widehat{\Sigma}}}_{ij} - u_{\widehat{\Sigma}} \brac{g_{\widehat{\Sigma}}}_{ij} -\frac{u_{\widehat{\Sigma}} \sscalar{p,\partial_i p} \sscalar{p, \partial_j p}}{1-\abs{p}^2} }, \label{GibS-2FF}
\end{align}
where the left-hand sides are evaluated at $P = \brac{p, p_{n+2}} \in \Sigma$ and the right-hand sides are evaluated at $p \in \widehat{\Sigma}$.
\end{lem}

\begin{proof}
The formula \eqref{GibS-lambda'} follows from the model \eqref{ModelHemisphere} of $\S^{n+1}_+$ in $\R^{n+2}$. For clarity, we view $p_{n+2}$ as the function $(1- \abs{p}^2)^{\frac{1}{2}}$ in the following calculation.
Differentiating $P = \brac{p, p_{n+2}}$ on $\Sigma$ and $\widehat{\Sigma}$ yields
\begin{equation}\label{fml:dP}
\partial_i P = \brac{\partial_i p, - \brac{1-\abs{p}^2}^{-\frac{1}{2}} \sscalar{p, \partial_i p}},
\end{equation}
and thus
\begin{align*}
\brac{g_{\Sigma}}_{ij} = \sscalar{\partial_i P, \partial_j P} = \brac{g_{\widehat{\Sigma}}}_{ij}+ \frac{\sscalar{p, \partial_i p} \sscalar{p, \partial_j p} }{1-\abs{p}^2}.
\end{align*}

Define $\nu_{\Sigma}$ as the right-hand side of \eqref{GibS-nu}. Note that $\sscalar{\nu_{\Sigma}, P} =0$, $\sscalar{\nu_{\Sigma}, \partial_i P} = 0$, and $\abs{\nu_{\Sigma}} =1$. This means that $\nu_{\Sigma}$ is a unit normal to $\Sigma$ at $P$. Moreover, since $\sscalar{d\pi \brac{\nu_{\Sigma}}, \nu_{\widehat{\Sigma}}} = \sqrt{1- u_{\widehat{\Sigma}}^2} >0$, we have that $\nu_{\Sigma}$ is outward pointing. Hence, the vector field $\nu_{\Sigma}$ given by \eqref{GibS-nu} is the outward unit normal to $\Sigma$.

Using $\lambda''+\lambda=0$, \eqref{GibS-lambda'} and \eqref{GibS-nu}, we have
\begin{equation*}
u_{\Sigma} = \sscalar{\lambda\brac{r} \partial_r, \nu_{\Sigma}} =  -\sscalar{\nabla^\S \lambda'\brac{r}, \nu_{\Sigma}} = -\sscalar{\brac{0,1}, \nu_{\Sigma}} = \brac{\frac{1-\abs{p}^2 }{1-u_{\widehat{\Sigma}}^2}}^{\frac{1}{2}} u_{\widehat{\Sigma}}.
\end{equation*}

Differentiating both sides of \eqref{fml:dP} yields
\begin{equation*}
\partial_j \partial_i P = \brac{\partial_j \partial_i p, -p_{n+2}^{-1} \sscalar{p, \partial_j \partial_i p} -p_{n+2}^{-1} \sscalar{\partial_i p, \partial_j p} - p_{n+2}^{-3} \sscalar{p, \partial_i p} \sscalar{p, \partial_j p}  }.
\end{equation*}
Together with \eqref{GibS-nu}, this gives
\begin{align*}
\brac{\II_{\Sigma}}_{ij} &= - \sscalar{\partial_j \partial_i P, \nu_{\Sigma}} \\
&= \brac{1-u_{\widehat{\Sigma}}^2}^{-\frac{1}{2}} \brac{\brac{\II_{\widehat{\Sigma}}}_{ij} - u_{\widehat{\Sigma}} \sscalar{\partial_i p, \partial_j p} - p_{n+2}^{-2} u_{\widehat{\Sigma}} \sscalar{p,\partial_i p} \sscalar{p, \partial_j p}}\\
&=\brac{1-u_{\widehat{\Sigma}}^2}^{-\frac{1}{2}} \brac{\brac{\II_{\widehat{\Sigma}}}_{ij} - u_{\widehat{\Sigma}} \brac{g_{\widehat{\Sigma}}}_{ij} -\frac{u_{\widehat{\Sigma}} \sscalar{p,\partial_i p} \sscalar{p, \partial_j p}}{1-\abs{p}^2} }.
\end{align*}
Thus, we obtain \eqref{GibS-2FF}. This completes the proof.
\end{proof}

The following proposition is the counterpart of Proposition \ref{prop:GibH-StaticC-Convex} for $\M^{n+1} = \S_+^{n+1}$.
\begin{prop}\label{prop:GibS-StaticC-Convex}
Let $\Sigma$ be a $C^2$ hypersurface in $\S_+^{n+1}$. Then $\Sigma$ is (strictly) static convex in $\S_+^{n+1}$ if and only if $\widehat{\Sigma}$ is (strictly) convex in $\R^{n+1}$.
\end{prop}
\begin{proof}
Taking $c=1$ in \eqref{def:S-Sigma} and using Lemma \ref{lem:GibS-GeoQuantitu}, we have
\begin{equation*}
\lambda'\brac{r} \brac{S_{\Sigma}}_{ij} = 
\lambda'\brac{r} \brac{\II_{\Sigma}}_{ij}+ u_{\Sigma} \brac{g_{\Sigma}}_{ij} = \brac{\frac{1-\abs{p}^2}{1-u_{\widehat{\Sigma}}^2}}^{\frac{1}{2}} \brac{\II_{\widehat{\Sigma}}}_{ij}.
\end{equation*}
Since $\lambda'\brac{r}>0$, the assertion follows.
\end{proof}

The following lemma is the analog of Lemma \ref{lem:Convex-Tensor-H} in $\S^{n+1}_+$.
\begin{lem}\label{lem:Convex-Tensor-S}
Let $\Omega$ be a bounded domain in $\S^{n+1}_+$ with $C^2$ boundary $\Sigma$, and let $K_{\Omega} = \pi \brac{\Omega}$. Assume that $\Sigma$ is strictly static convex and that $\II_{\Sigma} >\varepsilon g_{\Sigma}$ for some $\varepsilon \in \R$. Then the tensor
\begin{equation*}
\brac{h_{K_{\Omega}} \brac{1-h_{K_{\Omega}}^2}^{-\frac{1}{2}}+ \varepsilon} D^2 h_{K_{\Omega}} + \brac{1-h_{K_{\Omega}}^2}^{-\frac{3}{2}} Dh_{K_{\Omega}} \otimes D h_{K_{\Omega}}- \brac{\brac{1-h_{K_{\Omega}}^2}^{\frac{1}{2}}- \varepsilon h_{K_{\Omega}} } g_{\S^n}
\end{equation*}
is negative definite on $\S^n$.
\end{lem}

\begin{proof}
Using \eqref{GibS-2FF} and \eqref{GibS-metric}, we have
\begin{equation}\label{eq:II-metric-GibS}
\brac{\II_{\Sigma}- \varepsilon g_{\Sigma}}_{ij} = \frac{1}{\sqrt{1-u_{\widehat{\Sigma}}^2}} \brac{\II_{\widehat{\Sigma}}}_{ij} - \brac{\frac{u_{\widehat{\Sigma}}}{\sqrt{1- u_{\widehat{\Sigma}}^2 }}+ \varepsilon  }\brac{\brac{g_{\widehat{\Sigma}}}_{ij} + \frac{\sscalar{p, \partial_i p} \sscalar{p, \partial_j p} }{1-\abs{p}^2} }.
\end{equation}
By the assumption $\II_{\Sigma} >\varepsilon g_{\Sigma}$, the right-hand side of \eqref{eq:II-metric-GibS} is positive definite on $\widehat{\Sigma} = \partial K_{\Omega}$.

Write $h$ for $h_{K_{\Omega}}$. By Proposition \ref{prop:GibS-StaticC-Convex} and $K_{\Omega} \subset B^{n+1}$, we have $D^2 h + h g_{\S^n}>0$ and $\abs{Dh}^2+h^2<1$ on $\S^n$. If $h\brac{1-h^2}^{-\frac{1}{2}}+ \varepsilon \leq 0$ at $\theta \in \S^n$, then
\begin{align*}
&\brac{h \brac{1-h^2}^{-\frac{1}{2}} + \varepsilon} D^2 h+ \brac{1-h^2}^{-\frac{3}{2}} Dh \otimes Dh- \brac{\brac{1-h^2}^{\frac{1}{2}}-\varepsilon h } g_{\S^n} \\
&\leq -h\brac{h \brac{1-h^2}^{-\frac{1}{2}} + \varepsilon}g_{\S^n} + \brac{1-h^2}^{-\frac{3}{2}} Dh \otimes Dh- \brac{\brac{1-h^2}^{\frac{1}{2}}-\varepsilon h } g_{\S^n}\\
&= -\brac{1-h^2}^{-\frac{3}{2}} \brac{\brac{1-h^2} g_{\S^n}- Dh \otimes Dh} <0.
\end{align*} 
Hence, we can assume $h\brac{1-h^2}^{-\frac{1}{2}}+ \varepsilon >0$.
It follows from \eqref{eq:II-metric-GibS}, \eqref{eq:Rel-Supports} and \eqref{fml:2FF-D2h} that
\begin{equation*}
\brac{D^2 h + h \mathrm{Id}}^{-1} > \brac{\varepsilon \brac{1-h^2}^{\frac{1}{2}} + h} \brac{\mathrm{Id}+ \frac{Dh \otimes Dh}{1-\abs{Dh}^2-h^2} }.
\end{equation*}
Taking the inverse of both sides yields
\begin{equation*}
D^2 h+ h \mathrm{Id} <\frac{1}{\varepsilon \brac{1-h^2}^{\frac{1}{2}}+h} \brac{\mathrm{Id} - \frac{Dh \otimes Dh}{1-h^2}}.
\end{equation*}
Rearranging this inequality gives
\begin{equation*}
\brac{\varepsilon \brac{1-h^2}^{\frac{1}{2}}+h} D^2 h+ \brac{\varepsilon h\brac{1-h^2}^{\frac{1}{2}}+h^2-1} \mathrm{Id}+ \frac{Dh \otimes Dh}{1-h^2} <0.
\end{equation*}
Multiplying both sides by $\brac{1-h^2}^{-\frac{1}{2}}$, we obtain
the assertion.
\end{proof}

A proof similar to that of Lemma \ref{lem:regularity-Lomega-H} yields the following lemma. 
\begin{lem}\label{lem:regularity-Lomega-S}
Let $\Omega$ be a bounded domain in $\S^{n+1}_+$ with $C^3$ boundary $\Sigma$. Assume that $\II_{\Sigma} > \varepsilon g_{\Sigma}$ for some $\varepsilon \in \R$. Then the map $F : \Sigma \to \R^{n+1}$ defined by $F (P) = \pi \brac{\nu \brac{P} -\varepsilon P}$ is a $C^2$ immersion.
\end{lem}

\subsection{Projection of geodesic balls}
By Proposition \ref{prop:BallStaticConvex}, geodesic balls in $\M^{n+1}$ are strictly static convex. Then Proposition \ref{prop:GibH-StaticC-Convex} and Proposition \ref{prop:GibS-StaticC-Convex} imply that their images under the orthogonal projection $\pi$ are convex bodies with strictly convex boundaries in $\R^{n+1}$. In the following Lemma \ref{lem:GeoBall-Ellip-H} and Lemma \ref{lem:GeoBall-Ellip-S}, we establish the stronger result that these images are in fact ellipsoids.

\begin{lem}\label{lem:GeoBall-Ellip-H}
Let $\Omega \subset \H^{n+1}$ be a geodesic ball of radius $\rho$ centered at $X = \brac{x, x_{n+2}}$. Then $K_{\Omega} :=\pi\brac{\Omega}$ is an ellipsoid in $\R^{n+1}$. Moreover, if $X \neq \mathcal{O}$, then the ellipsoid $K_{\Omega}$ is centered at $\lambda'(\rho)x$, with one principal axis in the direction of $x$ of semi-axis length $\lambda(\rho)x_{n+2}$, and the remaining $n$ principal axes are orthogonal to $x$, each with a semi-axis length of $\lambda(\rho)$. If $X = \mathcal{O}$, then $\pi \brac{\Omega}$ is a centered Euclidean ball of radius $\lambda \brac{\rho}$. In either case, the support function of $K_{\Omega}$ is
\begin{equation}\label{fml:supp-Ellip-H}
h_{K_{\Omega}} \brac{\theta} = \lambda\brac{\rho} \sqrt{1+ \sscalar{\theta,x}^2} + \lambda'\brac{\rho} \sscalar{\theta,x}, \quad \theta \in \S^n.
\end{equation}
\end{lem}

\begin{proof}
If $X = \mathcal{O}$, then $K_{\Omega}$ is the Euclidean ball of radius $\lambda \brac{\rho}$ centered at the origin, since a point $(r, \theta)$ in the polar coordinate system  of $\H^{n+1}$ is identified with $(\lambda\brac{r}\theta, \lambda'\brac{r})$ on the hyperboloid.

Assume that $X \neq \mathcal{O}$, and let $y \in \pi \brac{\Omega}$. Then  $Y = \brac{y, y_{n+2}} \in \Omega$, where $y_{n+2} = \sqrt{1+\abs{y}^2}$. In view of \eqref{fml:GeoDis-H}, we have
\begin{equation}\label{condH-1}
\lambda'\brac{\rho} \geq -\sscalar{Y, X}_1=x_{n+2} y_{n+2} - \sscalar{x, y}.
\end{equation}
Write $y= \lambda'\brac{\rho} x + z + w$, where $z \in \R^{n+1}$ is parallel to $x$, and $w \in \R^{n+1}$ is perpendicular to $x$. Substituting this decomposition into \eqref{condH-1}, we obtain
\begin{align*}
x_{n+2} \brac{1+\abs{\lambda' \brac{\rho} x+ z}^2+\abs{w}^2 }^{\frac{1}{2}} = x_{n+2} y_{n+2} \leq \lambda'\brac{\rho} + \lambda'\brac{\rho} \abs{x}^2+ \sscalar{x, z}.
\end{align*}
Squaring both sides and simplifying yields
\begin{equation}\label{eq:Ellip-H}
x_{n+2}^{-2} \abs{z}^2+ \abs{w}^2 \leq \lambda \brac{\rho}^2.
\end{equation}
Thus, $K_{\Omega} = \pi \brac{\Omega}$ is the ellipsoid described in the statement.

Next, we calculate the support function $h_{K_{\Omega}}$. Let $\theta \in \S^n$. Let $\theta_1$ be the projection of $\theta$ on $\R x$, i.e., $\theta_1 = \sscalar{\theta, x} x/ \abs{x}^2$, and let $\theta_2 = \theta -\theta_1$. Then for any $y \in K_{\Omega}$, we have by using the Cauchy--Schwarz inequality and \eqref{eq:Ellip-H},
\begin{align*}
\sscalar{\theta, y- \lambda'\brac{\rho} x} &= \sscalar{\theta_1, z}+ \sscalar{\theta_2, w} \\
&\leq \brac{x_{n+2}^{-2} \abs{z}^2 + \abs{w}^2}^{\frac{1}{2}} \brac{x_{n+2}^2 \abs{\theta_1}^2+ \abs{\theta_2}^2 }^{\frac{1}{2}}\\
&\leq \lambda \brac{\rho} \brac{x_{n+2}^2 \abs{\theta_1}^2+ \abs{\theta_2}^2 }^{\frac{1}{2}}\\
&= \lambda\brac{\rho} \sqrt{1+ \sscalar{\theta, x}^2}.
\end{align*}
It is clear that the inequality above can attain equality by choosing suitable $y \in K_{\Omega}$. 

Consequently,
\begin{equation*}
h_{K_{\Omega}} \brac{\theta} = \sup_{y \in K_{\Omega}} \sscalar{\theta, y} = \lambda \brac{\rho} \sqrt{1+ \sscalar{\theta, x}^2} + \lambda'\brac{\rho} \sscalar{\theta, x}.
\end{equation*}
The assertion follows.
\end{proof}

The following Lemma \ref{lem:GeoBall-Ellip-S} is the counterpart of Lemma \ref{lem:GeoBall-Ellip-H} in $\S^{n+1}_+$.
\begin{lem}\label{lem:GeoBall-Ellip-S}
Let $\Omega \subset \S_+^{n+1}$ be a geodesic ball of radius $\rho$ centered at $P = \brac{p, p_{n+2}}$. Then $K_{\Omega} :=\pi\brac{\Omega}$ is an ellipsoid in $B^{n+1}$. Moreover, if $P \neq \mathcal{O}$, then the ellipsoid $K_{\Omega}$ is centered at $\lambda'(\rho)p$, with one principal axis in the direction of $p$ of semi-axis length $\lambda(\rho)p_{n+2}$, and the remaining $n$ principal axes are orthogonal to $p$, each with a semi-axis length of $\lambda(\rho)$. If $P = \mathcal{O}$, then $\pi \brac{\Omega}$ is a centered Euclidean ball of radius $\lambda \brac{\rho}$. In either case, the support function of $K_{\Omega}$ is
\begin{equation}\label{fml:supp-Ellip-S}
h_{K_{\Omega}} \brac{\theta} = \lambda\brac{\rho} \sqrt{1- \sscalar{\theta,p}^2} + \lambda'\brac{\rho} \sscalar{\theta,p}, \quad \theta \in \S^n.
\end{equation}
\end{lem}

\begin{proof}
If $P = \mathcal{O}$, then $K_{\Omega}$ is the Euclidean ball of radius $\lambda\brac{\rho}$ centered at the origin, since a point $(r, \theta)$ in the polar coordinate system of $\S^{n+1}_+$ is identified with $\brac{\lambda\brac{r} \theta, \lambda'\brac{r}}$ on the upper hemisphere.

Assume that $P \neq \mathcal{O}$. Let $q \in K_{ \Omega}$ and set $Q = (q, q_{n+2}) = \pi^{-1} \brac{q} $. Then, 
\begin{equation}\label{condS-1}
\lambda' \brac{\rho} \leq \sscalar{P, Q } = p_{n+2} q_{n+2} + \sscalar{p,q}.	
\end{equation}
Since $\Omega \subset \S^{n+1}_+$, we have $p_{n+2} > \lambda\brac{\rho}$, which implies $\abs{p} < \lambda'\brac{\rho}$. Write $q = \lambda'\brac{\rho} p + z+w$, where $z$ is parallel to $p$, and $w$ is perpendicular to $p$. Substituting this into \eqref{condS-1} and using $\abs{p}<\lambda'\brac{\rho}$, we obtain
\begin{equation*}
p_{n+2}^{-2} \abs{z}^2+ \abs{w}^2 \leq \lambda\brac{\rho}^2.
\end{equation*}
This yields the structure of the ellipsoid $K_{\Omega}$ described in the statement, and \eqref{fml:supp-Ellip-S} follows by a similar argument as that in Lemma \ref{lem:GeoBall-Ellip-H}. 
\end{proof}

\section{Steiner formula for the weighted volume} \label{sec:Steiner}
Let $\Omega$ be a bounded domain in $\M^{n+1}$. Let $\Omega_t$ denote the outer parallel set of $\Omega$ at distance $t \geq 0$, that is, $\Omega_t := \{X \in \M^{n+1};~\mathrm{dist}_{g_\M} \brac{X, \Omega} \leq t \}$. In this section, we derive the following Steiner formula \eqref{Fml:Steiner-WVol-CSC} for the weighted volume. The case $\varepsilon=1$ can be found in \cite{LX}*{Thm.~6.4}.

\begin{prop}\label{prop:SteinerFormulaWeightedCurvature}
Let $\Omega$ be a bounded domain in $\M^{n+1} \in \set{\H^{n+1}, \S^{n+1}_+}$ with $C^2$ boundary $\Sigma$. Let $\varepsilon \in \R$ be a constant. Then there exists a constant $\delta>0$ such that for any $t \in [0, \delta)$,
\begin{equation}\label{Fml:Steiner-WVol-CSC}
W_0^{\lambda'} \brac{\Omega_t} = \sum_{k=0}^{n+1} \binom{n+1}{k}  \brac{\lambda' \brac{t}+ \varepsilon \lambda\brac{t}}^{n+1-k} \lambda \brac{t}^k W_k^{\lambda', \varepsilon} \brac{\Omega}.
\end{equation}    
\end{prop}

\begin{proof}
For sufficiently small $t \geq 0$, $\Sigma_t =\partial \Omega_t$ is a $C^2$ hypersurface and evolves with velocity $\partial_t X = \nu$. Denote by $u_t: \Sigma_t \to \R$ the support function of $\Sigma_t$, and write $u$ for $u_{\Sigma}$. The evolution equations \eqref{fml:evl-Lam} and \eqref{fml:evl-U} then imply
\begin{equation}\label{fml:utStein}
u_t \brac{X_t} = \lambda \brac{t} \lambda'\brac{r} + \lambda' \brac{t} u = \lambda \brac{t} \brac{\lambda'\brac{r} - \varepsilon u} + \brac{\lambda'\brac{t} + \varepsilon \lambda \brac{t}} u,
\end{equation}
where $X_t$ denotes the position vector of $\Sigma_t = \partial \Omega_t$. From the geodesic equations \eqref{eq:geoH} and \eqref{eq:geoS}, we have
\begin{equation}\label{fml:PosNeigh}
X_t= \lambda' \brac{t} X + \lambda \brac{t} \nu_{\Sigma}, \quad X \in \Sigma.
\end{equation}
Note that $X_t \in \R^{n+1,1}$ when $\M^{n+1} = \H^{n+1}$, and $X_t \in \R^{n+2}$ when $\M^{n+1} = \S^{n+1}_+$. Moreover, by differentiating \eqref{fml:PosNeigh} on $\Sigma$ and taking the determinant in a local orthonormal frame, we have 
\begin{equation} \label{eq:pullbackVolForm}
X_t^* \brac{d\mu_{\Sigma_t}} = \det \brac{\lambda' \brac{t} g_{\Sigma} + \lambda \brac{t} \II_{\Sigma} } \, d\mu_{\Sigma},
\end{equation}
where we can write
\begin{equation*}
\lambda ' \brac{t} g_{\Sigma} + \lambda \brac{t} \II_{\Sigma}
= \brac{\lambda' \brac{t} +\varepsilon \lambda \brac{t}} g_{\Sigma} + \lambda \brac{t} \brac{\II_{\Sigma} - \varepsilon g_{\Sigma}}.
\end{equation*}
If $\M^{n+1} = \H^{n+1}$, we choose $\delta>0$ sufficiently small so that the outward normal exponential map of $\Sigma$ parametrizes $\Sigma_t$ for every $t \in [0, \delta)$; if $\M^{n+1} = \S^{n+1}_+$, we further require that $\delta$ is less than the distance from $\Omega$ to the boundary of $\S^{n+1}_+$.
Using \eqref{fml:utStein}, \eqref{eq:pullbackVolForm} and \eqref{Wk-EquiDef}, we obtain
\begin{align*}
W_0^{\lambda'} \brac{\Omega_t} &= \int_{\Sigma_t} u_t \, d\mu_{\Sigma_t} \\
&= \int_{\Sigma} \brac{\brac{\lambda' \brac{r}-\varepsilon u} \lambda \brac{t} +u \brac{\lambda' \brac{t}+\varepsilon \lambda \brac{t}} } \det \brac{\lambda'\brac{t} g_{\Sigma} + \lambda \brac{t} \II_{\Sigma}}d\mu_{\Sigma}\\
&= \sum_{k=0}^n \brac{\int_{\Sigma}  \brac{\lambda'\brac{r} -\varepsilon u} p_k \brac{\kappa-\varepsilon} \,d\mu_{\Sigma}} \binom{n}{k} \lambda\brac{t}^{k+1} \brac{\lambda' \brac{t}+\varepsilon \lambda \brac{t}}^{n-k}\\
&\qquad +\sum_{k=0}^n \brac{\int_{\Sigma} u p_k \brac{\kappa -\varepsilon} \, d\mu_{\Sigma}} \binom{n}{k} \lambda\brac{t}^k\brac{\lambda' \brac{t}+\varepsilon \lambda \brac{t}}^{n-k+1}\\
&= W_0^{\lambda', \varepsilon}\brac{\Omega} \brac{\lambda'\brac{t}+\varepsilon \lambda \brac{t}}^{n+1} + W_{n+1}^{\lambda',\varepsilon}\brac{\Omega} \lambda\brac{t}^{n+1} \\
&\qquad + \sum_{k=1}^n W_k^{\lambda', \varepsilon} \brac{\Omega} \brac{\binom{n}{k-1} +\binom{n}{k}}\lambda\brac{t}^k \brac{\lambda' \brac{t}+\varepsilon \lambda \brac{t}}^{n-k+1} \\
&= \sum_{k=0}^{n+1} \binom{n+1}{k} \brac{\lambda' \brac{t}+\varepsilon \lambda \brac{t}}^{n-k+1}\lambda \brac{t}^k W_k^{\lambda',\varepsilon} \brac{\Omega}.
\end{align*}
This completes the proof.
\end{proof}

\begin{rem}\label{rem:LinearIndep}
For $t \in (0, \delta)$,
$$\brac{\lambda' \brac{t}+ \varepsilon \lambda\brac{t}}^{n-k+1} \lambda \brac{t}^k = \lambda \brac{t}^{n+1} \brac{\lambda' \brac{t}/ \lambda \brac{t} + \varepsilon}^{n-k+1}.$$
Since $\lambda' \brac{t}/ \lambda \brac{t} + \varepsilon$ is not constant on $\brac{0, \delta}$, the functions $\brac{\lambda' \brac{t}+ \varepsilon \lambda\brac{t}}^{n-k+1} \lambda \brac{t}^k$ with $ k= 0, \ldots, n+1$ are linearly independent. 
\end{rem}
\section{Decomposition of projections of outer parallel sets} \label{sec:Decomp}
Let $\Omega$ be a bounded domain in $\H^{n+1}$ with boundary $\Sigma$. Define $K_{\Omega}$ as the domain enclosed by $\pi \brac{\Sigma}$ in $\R^{n+1} \times \set{0}$, that is, $K_{\Omega} = \pi \brac{\Omega}$. Similarly, when it is well defined, let $L_{\Omega}$ denote the domain enclosed by $\pi \brac{\nu \brac{\Sigma}}$. In view of the formula \eqref{fml:PosNeigh}
\begin{equation*}
X_t = \lambda'\brac{t} X+ \lambda \brac{t} \nu_{\Sigma}
\end{equation*}
for the position vector of $\Sigma_t= \partial \Omega_t$, one can expect that 
\begin{equation}
\pi \brac{\Omega_t} := K_{\Omega_t} = \lambda'\brac{t} K_{\Omega} + \lambda\brac{t} L_{\Omega} \label{sum-hope}
\end{equation}
when $K_{\Omega}$ and $L_{\Omega}$ are convex bodies in $\R^{n+1}$. By Proposition \ref{prop:GibH-StaticC-Convex}, the convexity of $K_{\Omega} \subset \R^{n+1}$ is guaranteed by the static convexity of $\Omega \subset \H^{n+1}$. In Lemma \ref{lem:DecompKt-H} below, we show that $L_{\Omega}=L_{\Omega,0}$ is a convex body when $\Omega$ is convex and static convex. The decomposition \eqref{sum-hope} will be used to prove the case $\varepsilon=0$ of Theorem \ref{MainThm-1} later. However, this case recovers the second Minkowski type inequality \eqref{Ineq:Xia16} only under the additional assumption that $\Omega$ is convex.

To remove this additional convexity assumption, we introduce an alternative decomposition of $K_{\Omega_t}$ using the vector field $\nu_{\Sigma}+X$. This is a future-directed null normal of $\Sigma \subset \R^{n+1,1}$ in the context of general relativity; see \cite{BW14}. Define $T_{\Omega}$ as the domain enclosed by $\pi \brac{\brac{\nu_{\Sigma}+X} \brac{\Sigma}}$. As shown in Lemma \ref{lem:DecompKt-H} below, $T_{\Omega} = L_{\Omega,-1}$ is well defined and convex whenever $\Omega$ is static convex. We then express $K_{\Omega_t}$ as a Minkowski combination of $K_{\Omega}$ and $T_{\Omega}$ and will use this decomposition to prove the case $\varepsilon = -1$ of Theorem \ref{MainThm-1}. Together with Corollary \ref{cor-LWXObser}, this recovers the inequality \eqref{Ineq:Xia16} under the assumption that $\Omega$ is static convex.

To unify these two decompositions of $K_{\Omega_t}$, we introduce a parameter $\varepsilon$ in the following Lemma \ref{lem:DecompKt-H}. We also study the corresponding decomposition of $K_{\Omega_t}$ for $\Omega \subset \S^{n+1}_+$ in Lemma \ref{lem:DecompKt-S} below.

\begin{lem} \label{lem:DecompKt-H}
Let $\Omega$ be a bounded static convex domain in $\H^{n+1}$ with $C^3$ boundary $\Sigma$. Assume that $\II_{\Sigma} \geq \varepsilon g_{\Sigma}$ for some $\varepsilon \in \R$. Then there exists $\delta>0$ such that $K_{\Omega_t} := \pi \brac{\Omega_t}$ is a convex body in $\R^{n+1}$ with support function $\lambda'\brac{t} h_{K_{\Omega}}+ \lambda \brac{t} \sqrt{1+ h_{K_\Omega}^2 }$ for any $t \in [0,\delta)$. Moreover, we have
\begin{equation}\label{decomp-KT}
K_{\Omega_t} = \brac{\lambda'\brac{t} + \varepsilon \lambda \brac{t}} K_{\Omega} + \lambda \brac{t} L_{\Omega, \varepsilon},
\end{equation}
where $L_{\Omega, \varepsilon}$ is a convex body in $\R^{n+1}$ with support function 
\begin{equation} \label{eq:SupportLeps-H}
h_{L_{\Omega, \varepsilon}} = \sqrt{1+h_{K_{\Omega}}^2}- \varepsilon h_{K_{\Omega}}.
\end{equation}
The boundary of $L_{\Omega, \varepsilon}$ is given by $\set{ \pi \brac{\nu_{\Sigma}\brac{X} - \varepsilon X}; X \in \Sigma}$.
\end{lem}

\begin{proof}
We first consider the case where $\Omega$ is strictly static convex with $\II_{\Sigma} > \varepsilon g_{\Sigma}$ on $\Sigma$. By Proposition \ref{prop:GibH-StaticC-Convex}, $K_{\Omega} = \pi \brac{\Omega}$ is a convex body with strictly convex boundary in $\R^{n+1}$. 

Let $X= \brac{x,x_{n+2}} \in \Sigma$, and let $\theta = \nu_{\partial K_{\Omega}} \brac{x} \in \S^n$. It follows from the strict convexity of $K_{\Omega}$ that
\begin{equation} \label{fml:pi(X)}
\pi \brac{X} = D h_{K_{\Omega}} + h_{K_{\Omega}} \theta.
\end{equation}
Then \eqref{eq:GibH-nu} and \eqref{eq:Rel-Supports} yield
\begin{equation}\label{fml:pi(nu)-Support}
\pi \brac{\nu_{\Sigma} \brac{X}} = \frac{\theta + h_{K_{\Omega}} \brac{D h_{K_{\Omega}}+h_{K_{\Omega}} \theta  } }{\brac{1+h_{K_{\Omega}}^2}^{\frac{1}{2}}}= D \brac{1+h_{K_{\Omega}}^2}^{\frac{1}{2}} +\brac{1+h_{K_{\Omega}}^2}^{\frac{1}{2}} \theta. 
\end{equation}
Let $f = \sqrt{1+h_{K_{\Omega}}^2} - \varepsilon h_{K_{\Omega}}$. The above formulas yield
\begin{equation}\label{eq:pi(nu-X)}
\pi \brac{\nu_{\Sigma} \brac{X} - \varepsilon X} = Df + f \theta
\end{equation}
and
\begin{align}
\pi \brac{X_t} &= \pi \brac{\lambda' \brac{t} X+ \lambda \brac{t} \nu_{\Sigma}} = \brac{\lambda' \brac{t}+ \varepsilon \lambda \brac{t}} \brac{D h_{K_\Omega}+ h_{K_{\Omega}} \theta} + \lambda \brac{t} \brac{D f + f \theta}.  \label{eq:pi(Xt)}
\end{align}
By a direct calculation,
\begin{align*}
D^2 f + f g_{\S^n} &=\frac{1}{\brac{1+h_{K_{\Omega}}^2}^{\frac{1}{2}} } \brac{h_{K_{\Omega}} D^2 h_{K_{\Omega}}+ \frac{Dh_{K_{\Omega}} \otimes Dh_{K_{\Omega}}}{1+h_{K_{\Omega}}^2}+ \brac{1+h_{K_{\Omega}}^2} g_{\S^n} }\\
&\qquad  - \varepsilon \brac{D^2 h_{K_{\Omega}} + h_{K_{\Omega}} g_{\S^n}}.
\end{align*}
It follows from Lemma \ref{lem:Convex-Tensor-H} that $D^2 f+ f g_{\S^n}$ is positive definite on $\S^n$. This implies that $f$ is the support function of a convex body $L_{\Omega, \varepsilon}$ with strictly convex boundary. Then $\eqref{eq:pi(nu-X)}$ shows that $L_{\Omega, \varepsilon}$ is exactly the domain enclosed by $\set{\pi \brac{\nu_{\Sigma} \brac{X} -\varepsilon X};~ X \in \Sigma}$. Let $\delta>0$ such that $\lambda'\brac{t}+\varepsilon \lambda\brac{t} \geq 0$ for any $t \in [0, \delta)$. The desired decomposition of $K_{\Omega_t}$ in \eqref{decomp-KT} therefore follows from \eqref{eq:pi(Xt)}.

For the general case where $\Omega$ is static convex with $\II_{\Sigma} \geq \varepsilon g_{\Sigma}$ on $\Sigma$, the assertion follows from an approximation argument. By Lemma \ref{lem:approximation-H}, there exists a sequence of strictly static convex domains $\Omega_i$ with boundary $\Sigma_i$ satisfying $\II_{\Sigma_i}> \varepsilon g_{\Sigma_i}$, and $\Sigma_i \to \Sigma$ in $C^2$-topology as $i \to \infty$.
In this case, $K_{(\Omega_i)_t}$, $K_{\Omega_i}$, and $L_{\Omega_i, \varepsilon}$ are convex bodies with strictly convex $C^2$ boundaries that satisfy the results in this lemma. 
As $i \to \infty$, $K_{\Omega_i}$ converges to $K_{\Omega}$, and $L_{\Omega_i, \varepsilon}$ converges to a convex body $L_{\Omega, \varepsilon}$, both in the Hausdorff metric. This
yields the decomposition \eqref{decomp-KT}. Moreover, the uniform convergence of the corresponding support functions on $\S^n$ implies that $h_{L_{\Omega, \varepsilon}}$ and $h_{K_{\Omega}}$ satisfy \eqref{eq:SupportLeps-H}.
\end{proof}

\begin{lem} \label{lem:DecompKt-S}
Let $\Omega$ be a bounded static convex domain in $\S^{n+1}_+$ with $C^3$ boundary $\Sigma$. Assume that $\II_{\Sigma} \geq \varepsilon g_{\Sigma}$ for some $\varepsilon \in \R$. Then there exists $\delta>0$ such that $K_{\Omega_t} := \pi \brac{\Omega_t}$ is a convex body in $B^{n+1}$ with support function $\lambda'\brac{t} h_{K_{\Omega}}+ \lambda \brac{t} \sqrt{1- h_{K_\Omega}^2}$ for any $t \in [0,\delta)$. Moreover, we have
\begin{equation}\label{decomp-KT-S}
K_{\Omega_t} = \brac{\lambda'\brac{t} + \varepsilon \lambda \brac{t}} K_{\Omega} + \lambda \brac{t} L_{\Omega, \varepsilon},
\end{equation}
where $L_{\Omega, \varepsilon}$ is a convex body in $B^{n+1}$ with support function 
\begin{equation}\label{eq:SupportLeps-S}
h_{L_{\Omega, \varepsilon}} = \sqrt{1-h_{K_{\Omega}}^2} - \varepsilon h_{K_{\Omega}}.
\end{equation}
The boundary of $L_{\Omega, \varepsilon}$ is given by $\set{ \pi \brac{\nu\brac{P} - \varepsilon P}; P \in \Sigma}$.
\end{lem}

\begin{proof}
Consider the case where $\Omega$ is strictly static convex with $\II_{\Sigma} > \varepsilon g_{\Sigma}$ on $\Sigma$. Let $P \in \Sigma$, and let $\theta = \nu_{\partial K_{\Omega}} \brac{p} \in \S^n$.
By Proposition \ref{prop:GibS-StaticC-Convex}, $K_{\Omega}$ is strictly convex. Then
\begin{equation*}
\pi \brac{P} = D h_{K_{\Omega}} +h_{K_{\Omega}} \theta.
\end{equation*}
Using \eqref{GibS-nu} and \eqref{eq:Rel-Supports}, we have
\begin{equation*}
\pi \brac{\nu_{\Sigma} \brac{P}} = \frac{\theta -h_{K_{\Omega}} \brac{D h_{K_{\Omega}}+h \theta} }{\brac{1-h_{K_{\Omega}}^2}^{\frac{1}{2}}} = D \brac{1-h_{K_{\Omega}}^2}^{\frac{1}{2}} +  \brac{1-h_{K_{\Omega}}^2}^{\frac{1}{2}} \theta.
\end{equation*}
Let $f= \sqrt{1-h_{K_{\Omega}}^2}- \varepsilon h_{K_{\Omega}}$. Then
\begin{equation*}
\pi \brac{\nu_{\Sigma} \brac{P}- \varepsilon P } = Df +f \theta
\end{equation*}
and
\begin{equation*}
\pi \brac{P_t} = \brac{\lambda' \brac{t} + \varepsilon \lambda \brac{t}} \brac{D h_{K_{\Omega}} + h_{K_{\Omega}} \theta} + \lambda \brac{t} \brac{Df + f \theta}.
\end{equation*}
A direct calculation yields
\begin{align*}
D^2 f+ f g_{\S^n} &= \frac{-1}{\brac{1-h_{K_{\Omega}}^2}^{\frac{1}{2}}} \brac{h_{K_{\Omega}} D^2 h_{K_{\Omega}} + \frac{D h_{K_{\Omega}} \otimes Dh_{K_{\Omega}} }{1-h_{K_{\Omega}}^2} - \brac{1-h_{K_{\Omega}}^2} g_{\S^n} }\\
&\qquad - \varepsilon \brac{D^2 h_{K_{\Omega}} +h_{K_{\Omega}} g_{\S^n}}.
\end{align*}
It follows from Lemma \ref{lem:Convex-Tensor-S} that $D^2 f+ f g_{\S^n}$ is positive definite on $\S^n$. Hence $\set{\pi \brac{\nu\brac{P}- \varepsilon P};~P \in \Sigma}$ encloses a convex body $L_{\Omega, \varepsilon}$ with support function $f$. Moreover, the decomposition \eqref{decomp-KT-S} holds for small $t \geq 0$. For the general case where $\Omega$ is static convex with $\II_{\Sigma} \geq \varepsilon g_{\Sigma}$ on $\Sigma$, the conclusion follows from an approximation argument, where the approximating hypersurfaces are given by Lemma \ref{lem:approximation-S}. 
\end{proof}

\section{Proof of Theorems \ref{MainThm-1} and \ref{Thm:FullRangeHuLi}} \label{sec:ProofThm}
\subsection{Proof of Theorem \ref{MainThm-1}}
In this subsection, we first prove Theorem \ref{MainThm-1} for the case $\M^{n+1} = \H^{n+1}$ and then for the case $\M^{n+1} = \S^{n+1}_+$.

\textbf{Case 1: $\M^{n+1} = \H^{n+1}$.}
Let $\Omega$ be a bounded static convex domain in $\H^{n+1}$ with $C^3$ boundary $\Sigma$. Let $\delta>0$ such that the Steiner type formula \eqref{Fml:Steiner-WVol-CSC} for the weighted volume $W_0^{\lambda'} \brac{\Omega_t}$ holds for any $t \in [0, \delta)$. In addition, we require $\lambda'\brac{t}+ \varepsilon \lambda \brac{t}>0$ for $t \in [0,\delta)$.
 
Since $\Omega$ is static convex with $\II_{\Sigma} \geq \varepsilon g_{\Sigma}$ on $\Sigma$, we can use Lemma \ref{lem:DecompKt-H} to decompose $K_{\Omega_t} = \pi \brac{\Omega_t}$ into linear Minkowski combinations of convex bodies as
\begin{equation}\label{eq:decomp-proof}
K_{\Omega_t} = \brac{\lambda' \brac{t}+\varepsilon \lambda \brac{t}} K_{\Omega} + \lambda \brac{t} L_{\Omega, \varepsilon},
\end{equation}
where $L_{\Omega,\varepsilon}$ is the convex body with boundary $\set{\pi \brac{\nu_{\Sigma}\brac{X}- \varepsilon X};~X \in \Sigma}$.

Using \eqref{eq:Vol-WeightedVol} and applying the polynomial expansion \eqref{fml:ClassicalSteiner} to the decomposition \eqref{eq:decomp-proof}, we have
\begin{align}
W_0^{\lambda'} \brac{\Omega_t} &= \brac{n+1}\Vol \brac{K_{\Omega_t}} \nonumber\\
&=\brac{n+1} \sum_{k=0}^{n+1} \binom{n+1}{k} V \brac{K_{\Omega}, n+1-k; L_{\Omega,\varepsilon}, k} \brac{\lambda' \brac{t}+\varepsilon \lambda \brac{t}}^{n+1-k} \lambda \brac{t}^k. \label{eq:VolKtExpan}
\end{align}
Recall from Remark \ref{rem:LinearIndep} that the functions $\brac{\lambda' \brac{t}+\varepsilon \lambda \brac{t}}^{n-k+1} \lambda \brac{t}^k$ are linearly independent on $[0, \delta)$ for $k=0, \ldots n+1$.
Comparing \eqref{eq:VolKtExpan} with the Steiner type formula for $W_0^{\lambda'} \brac{\Omega_t}$ in \eqref{Fml:Steiner-WVol-CSC}, we have
\begin{equation}\label{eq:rel-WeightCurvInteg-MixedVolume}
W_{k}^{\lambda',\varepsilon} \brac{\Omega} = \brac{n+1} V \brac{K_{\Omega}, n+1-k; L_{\Omega,\varepsilon}, k}, \quad  k=0, \ldots n+1.
\end{equation}
Since the mixed volume on the right-hand side of \eqref{eq:rel-WeightCurvInteg-MixedVolume} is nonnegative, we have $W_{k}^{\lambda', \varepsilon} \brac{\Omega} \geq 0$. Inserting \eqref{eq:rel-WeightCurvInteg-MixedVolume} into the Alexandrov--Fenchel inequality \eqref{ineq:AFReal}, we obtain
\begin{equation*}
W_{k}^{\lambda',\varepsilon} \brac{\Omega}^2 \geq W_{k-1}^{\lambda',\varepsilon} \brac{\Omega} W_{k+1}^{\lambda',\varepsilon} \brac{\Omega}, \quad k=1, \ldots, n.
\end{equation*}
Thus, the desired inequality \eqref{MainIneq:AlexFen} holds for the case $\M^{n+1} = \H^{n+1}$.

Now we classify the equality cases of the inequality \eqref{MainIneq:AlexFen}. If $\II_{\Sigma} > \varepsilon g_{\Sigma}$, then Lemma \ref{lem:regularity-Lomega-H} implies that $\partial L_{\Omega, \varepsilon} =\set{\pi \brac{\nu \brac{X} -\varepsilon X}; X \in \Sigma}$ is $C^2$ smooth. Thus $K_{\Omega}$ and $L_{\Omega, \varepsilon}$ are both convex bodies with $C^2$ boundaries and, in particular, have nonempty interiors. Assume that equality holds in \eqref{MainIneq:AlexFen}. Then by \eqref{eq:rel-WeightCurvInteg-MixedVolume} and Theorem \ref{thm:InSch14}, $K_{\Omega}$ and $L_{\Omega,\varepsilon}$ are homothetic. That is, there exist constant $\alpha >0$ and $v \in \R^{n+1}$ such that 
\begin{equation} \label{eq:homoKL}
L_{\Omega, \varepsilon} = \alpha K_{\Omega} + v.
\end{equation}
Using \eqref{eq:SupportLeps-H} and \eqref{eq:homoKL}, we have
\begin{equation}\label{eq:HomoSupp-H}
\sqrt{1+ h_{K_{\Omega}} \brac{\theta}^2} = \brac{\alpha+ \varepsilon}h_{K_{\Omega}} \brac{\theta} + \sscalar{v, \theta}, \quad \theta \in \S^n.
\end{equation}
The solutions to the above equation are
\begin{equation}\label{eq:SolCandi}
h_{K_{\Omega}} \brac{\theta} = \frac{-\brac{\alpha+\varepsilon} \sscalar{v, \theta} \pm \sqrt{\brac{\alpha+\varepsilon}^2-1+ \sscalar{v, \theta}^2 } }{\brac{\alpha+ \varepsilon}^2-1}.
\end{equation}
By taking $\theta_0 \in \S^n$ such that $\sscalar{v, \theta_0} \leq 0$ in \eqref{eq:HomoSupp-H}, we see that $\abs{h_{K_{\Omega}} \brac{\theta_0} } < \brac{\alpha+ \varepsilon} h_{K_{\Omega}} \brac{\theta_0}$, which implies $\abs{\alpha + \varepsilon} >1$. Since $K_{\Omega}$ has nonempty interior, we have $h_{K_{\Omega}} \brac{\theta} + h_{K_{\Omega}} \brac{-\theta}>0$. Then by \eqref{eq:HomoSupp-H}, we have $\alpha + \varepsilon >0$. Thus $\alpha+\varepsilon >1$, and there exists a unique $\rho >0$ such that $\alpha + \varepsilon = \lambda' \brac{\rho} / \lambda \brac{\rho}$. Using the fact $h_{K_{\Omega}} \brac{\theta} + h_{K_{\Omega}} \brac{-\theta}>0$ again to \eqref{eq:SolCandi}, we obtain
\begin{align*}
h_{K_{\Omega}} \brac{\theta} &= \frac{-\brac{\alpha+\varepsilon} \sscalar{v, \theta} + \sqrt{\brac{\alpha+\varepsilon}^2-1+ \sscalar{v, \theta}^2 } }{\brac{\alpha+ \varepsilon}^2-1}\\
&= \lambda'\brac{\rho} \sscalar{-\lambda\brac{\rho} v, \theta} +\lambda\brac{\rho} \sqrt{1+ \sscalar{-\lambda\brac{\rho} v, \theta}^2} .
\end{align*}
By Lemma \ref{lem:GeoBall-Ellip-H}, $\Omega \subset \H^{n+1}$ is a geodesic ball of radius $\rho$ centered at $\pi^{-1} \brac{- \lambda \brac{\rho} v}$.

Conversely, we let $\Omega$ be a geodesic ball in $\H^{n+1}$ centered at $X = \brac{x, x_{n+2}}$ with radius $\rho$ satisfying $\II_{\Sigma} > \varepsilon g_{\Sigma}$, equivalently, $\lambda'\brac{\rho} / \lambda \brac{\rho} > \varepsilon$. Then by \eqref{fml:supp-Ellip-H} and \eqref{eq:SupportLeps-H}, we have
\begin{align*}
h_{K_{\Omega}} \brac{\theta}&= \lambda\brac{\rho} \sqrt{1+\sscalar{\theta, x}^2} + \lambda'\brac{\rho} \sscalar{\theta, x}, \\
h_{L_{\Omega, \varepsilon}}\brac{\theta} &=\brac{\lambda'\brac{\rho}-\varepsilon \lambda \brac{\rho}}\sqrt{1+\sscalar{\theta, x}^2 } + \brac{\lambda\brac{\rho} -\varepsilon \lambda'\brac{\rho}}\sscalar{\theta, x}.
\end{align*}
Thus
\begin{equation*}
h_{L_{\Omega, \varepsilon}} = \brac{\frac{\lambda'\brac{\rho}}{\lambda\brac{\rho}} -\varepsilon} h_{K_{\Omega}} +\scalar{\theta, - \frac{x}{\lambda \brac{\rho}}},
\end{equation*}
which means that $L_{\Omega,\varepsilon}$ and $K_{\Omega}$ are homothetic. Due to the multi-linearity and translation invariance of the mixed volume, we have
\begin{equation*}
V \brac{K_{\Omega}; n+1-k; L_{\Omega,\varepsilon}, k}^2 = V \brac{K_{\Omega}, n+2-k; L_{\Omega, \varepsilon}, k-1} V \brac{K_{\Omega}, n-k; L_{\Omega, \varepsilon}, k+1}.
\end{equation*}
Then by \eqref{eq:rel-WeightCurvInteg-MixedVolume}, equality holds in \eqref{MainIneq:AlexFen}.

Consequently, when $\Omega$ is static convex with $\II_{\Sigma}> \varepsilon g_{\Sigma}$ on $\Sigma = \partial \Omega$, equality holds in \eqref{MainIneq:AlexFen} if and only if $\Omega \subset \H^{n+1}$ is a geodesic ball with radius $\rho$ satisfying $\lambda'\brac{\rho}/ \lambda \brac{\rho} >\varepsilon$. This completes the proof of Theorem \ref{MainThm-1} for the case $\M^{n+1} = \H^{n+1}$.

\textbf{Case 2: $\M^{n+1} = \S^{n+1}_+$.} For the proof of the inequality \eqref{MainIneq:AlexFen} in the case $\M^{n+1}= \S^{n+1}_+$, there is only one technical difference compared to the case $\M^{n+1} =\H^{n+1}$: instead of using Lemma \ref{lem:DecompKt-H}, we use Lemma \ref{lem:DecompKt-S} to obtain the decomposition \eqref{eq:decomp-proof}. 

For the classification of the equality cases of \eqref{MainIneq:AlexFen} when $\M^{n+1} = \S^{n+1}_+$, we provide a brief proof. Since $\II_{\Sigma}> \varepsilon g_{\Sigma}$, Lemma \ref{lem:regularity-Lomega-S} implies that $L_{\Omega, \varepsilon}$ is a convex body with $C^2$ boundary. Due to Theorem \ref{thm:InSch14}, if equality holds in \eqref{MainIneq:AlexFen}, then $L_{\Omega, \varepsilon} = \alpha K_{\Omega} +v$ for some $\alpha>0$ and $v \in \R^{n+1}$. Then we apply \eqref{eq:SupportLeps-S} to this equation and obtain
\begin{equation}\label{eq:HomoSupp-S}
\sqrt{1-h_{K_{\Omega}} \brac{\theta}^2} = \brac{\alpha+\varepsilon} h_{K_{\Omega}} \brac{\theta} +\sscalar{v, \theta}, \quad \theta \in \S^n.
\end{equation}
Since $K_{\Omega}$ has nonempty interior, it follows from $h_{K_{\Omega}} \brac{\theta} + h_{K_{\Omega}} \brac{-\theta} >0$ and the above \eqref{eq:HomoSupp-S} that $\alpha +\varepsilon >0$. Then there exists a unique $\rho \in (0, \pi/2)$ such that $\alpha+\varepsilon = \lambda'\brac{\rho} / \lambda \brac{\rho}$. Solving the equation \eqref{eq:HomoSupp-S} gives
\begin{equation}\label{eq:hK-HomoS}
h_{K_{\Omega}} \brac{\theta} = \lambda'\brac{\rho} \sscalar{p, \theta} + \lambda\brac{\rho} \sqrt{1-\sscalar{p, \theta}^2},
\end{equation}
where $p := - \lambda \brac{\rho} v$. 

We claim that $\abs{p} < \lambda' \brac{\rho}$. Suppose, to the contrary, that $\abs{p} \geq \lambda' \brac{\rho}$. Then there exists $\theta_0 \in \S^n$ such that $\sscalar{p, \theta_0} = \lambda' \brac{\rho}$. Inserting this into \eqref{eq:hK-HomoS}, we have $h_{K_{\Omega}} \brac{\theta_0} =1$, which contradicts the assumption $\Omega \subset \S^{n+1}_+$.

It follows from $\abs{p} <\lambda'\brac{\rho}$ that a geodesic ball in $\S^{n+1}$ of radius $\rho$ centered at $\pi^{-1} \brac{p}$ is contained in $\S^{n+1}_+$. In view of \eqref{eq:hK-HomoS} Lemma \ref{lem:GeoBall-Ellip-S}, we have that $\Omega$ is exactly the geodesic ball described above. Moreover, as $\alpha>0$, we have $\lambda'\brac{\rho}/ \lambda \brac{\rho} > \varepsilon$, which implies $\II_{\Sigma} > \varepsilon g_{\Sigma}$. Conversely, for any geodesic ball $\Omega$ in $\S^{n+1}_+$ of radius $\rho$ satisfying $\II_{\Sigma} > \varepsilon g_{\Sigma}$ on its boundary, one can use \eqref{fml:supp-Ellip-S} and \eqref{eq:SupportLeps-S} to verify that the corresponding $K_{\Omega}$ and $L_{\Omega, \varepsilon}$ are homothetic, and thus equality holds in \eqref{MainIneq:AlexFen}. This completes the proof of Theorem \ref{MainThm-1}.

\begin{rem}\label{rem:EqCase-StrictSC}
To provide geometric intuition, we present a proof of the equality case of \eqref{MainIneq:AlexFen} for $\Omega \subset \H^{n+1}$ under the additional assumption that $\Omega$ is strictly static convex, illustrating how the equation \eqref{eq:homoKL} implies that $\Omega$ is a geodesic ball. 

By the argument in the proof of Lemma \ref{lem:DecompKt-H}, $\partial K_{\Omega}$ and $\partial L_{\Omega, \varepsilon}$ are both strictly convex. In view of \eqref{fml:pi(X)} and \eqref{fml:pi(nu)-Support}, for any $\theta \in \S^n$, there exists a unique $X \in \Sigma$ such that the outward unit normal of $\partial K_{\Omega}$ at $\pi \brac{X}$ and that of $\partial L_{\Omega, \varepsilon}$ at $\pi \brac{\nu \brac{X} - \varepsilon X}$ are both $\theta$. This yields a pointwise version of \eqref{eq:homoKL}:
\begin{equation}\label{eq:PointwiseHomoKL}
\pi \brac{\nu \brac{X}-\varepsilon X} = \alpha \pi \brac{X}+ v, \quad  X \in \Sigma.
\end{equation}
Differentiating both sides of \eqref{eq:PointwiseHomoKL} yields
\begin{equation*}
\brac{\brac{\II_{\Sigma}}_i{}^j -(\alpha+\varepsilon) \delta_i{}^j} \pi \brac{\partial_j X} =0.
\end{equation*}
Since $\partial_j X \in T_X \H^{n+1}$ is spacelike, $\Sigma$ must be a geodesic sphere with all principal curvatures equal to $(\alpha+\varepsilon)$.
\end{rem}

\begin{proof}[Proof of Remark~\ref{Rmk:2}~\eqref{rmk2.2}]
If $\dim L_{\Omega, \varepsilon}=n+1$, then  \eqref{eq:rel-WeightCurvInteg-MixedVolume} implies that $W_l^{\lambda', \varepsilon} \brac{\Omega}>0$ for all $l=0, \ldots, n+1$. By \eqref{MainIneq:AlexFen}, $\log (W_l^{\lambda',\varepsilon} \brac{\Omega})$ is concave in $l$. This implies \eqref{ineq:GeneralV-AFmain}.

If $\dim L_{\Omega, \varepsilon}=d<n+1$, then Theorem \ref{thm:PosMixVol} implies that $V \brac{K_{\Omega},n+1-l; L_{\Omega,\varepsilon},l} >0$ precisely when $0 \leq l \leq d$. By \eqref{eq:rel-WeightCurvInteg-MixedVolume}, we have $W_l^{\lambda', \varepsilon} \brac{\Omega} >0$ precisely when $0 \leq l \leq d$. Hence $\log (W_l^{\lambda', \varepsilon} \brac{\Omega})$ is concave in $l$ for $0\leq l \leq d$, and $W_l^{\lambda', \varepsilon}\brac{\Omega} =0$ for $d+1 \leq l \leq n+1$. Thus \eqref{ineq:GeneralV-AFmain} follows from the concavity when $k \leq d$ and is trivial for $k>d$.
\end{proof}

\subsection{Proof of Theorem \ref{Thm:FullRangeHuLi}}
By Lemma \ref{lem:DecompKt-H}, $L_{\Omega,\varepsilon}$ is a convex body in $\R^{n+1}$ with $h_{L_{\Omega,\varepsilon}} = \sqrt{1+h_{K_{\Omega}}^2} - \varepsilon h_{K_{\Omega}}$. Since $\varepsilon \leq 0$, we have for any $\theta \in \S^n$,
\begin{equation*}
h_{L_{\Omega,\varepsilon}} \brac{\theta} + h_{L_{\Omega,\varepsilon}} \brac{-\theta} \geq 2 -\varepsilon \brac{h_{K_{\Omega}} \brac{\theta}+ h_{K_{\Omega}}\brac{-\theta}} \geq 2.
\end{equation*}
This means that the width of $L_{\Omega,\varepsilon}$ in each direction is at least $2$, and in particular, the interior of $L_{\Omega,\varepsilon}$ is nonempty. By \eqref{eq:rel-WeightCurvInteg-MixedVolume}, we have $W_i^{\lambda',\varepsilon} \brac{\Omega}>0$ for all $i=0, \ldots, n+1$.

The remainder of the proof of Theorem \ref{Thm:FullRangeHuLi} is divided into three steps.

\textbf{Step 1.} We start by establishing the inequality \eqref{MainIneq:ShifHuLi-H} in the case that $k=n+1$, $l=n$.  The identity \eqref{eq:rel-WeightCurvInteg-MixedVolume} yields
\begin{equation*}
W_{n+1}^{\lambda',\varepsilon} \brac{\Omega} =\brac{n+1} \Vol \brac{L_{\Omega,\varepsilon}}>0, \quad W_n^{\lambda',\varepsilon}  \brac{\Omega} =\brac{n+1} V \brac{K_{\Omega}, 1; L_{\Omega,\varepsilon}, n}>0. 
\end{equation*}
Then, using  \eqref{eq:VolMeasure} and \eqref{eq:AreaMeasure}, we have
\begin{align}
W_{n+1}^{\lambda',\varepsilon} \brac{\Omega} &= \int_{\S^n} h_{L_{\Omega,\varepsilon}} \, dS_{L_{\Omega,\varepsilon}} = \int_{\S^n} \brac{\brac{1+h_{K_{\Omega}}^2}^{\frac{1}{2}}-\varepsilon h_{K_{\Omega}}} \, dS_{L_{\Omega,\varepsilon}}, \label{eq:Wn+1SuppK}\\
W_n^{\lambda',\varepsilon} \brac{\Omega} &= \int_{\S^n} h_{K_{\Omega}} \, dS_{L_{\Omega,\varepsilon}}, \label{eq:WnSuppK}
\end{align}
where $dS_{\Omega,\varepsilon}$ is the surface area measure of $L_{\Omega,\varepsilon}$.
Since $f \brac{s} = \brac{1+s^2}^{\frac{1}{2}}$ is strictly convex on $\R$, applying Jensen's inequality to the right-hand side of \eqref{eq:Wn+1SuppK} and using \eqref{eq:WnSuppK}, we obtain
\begin{align}
W_{n+1}^{\lambda',\varepsilon} \brac{\Omega} &\geq \brac{\brac{\int_{\S^n} h_{K_{\Omega}} \, dS_{L_{\Omega,\varepsilon}}}^2+ \mathrm{Area} (\partial L_{\Omega,\varepsilon})^2}^{\frac{1}{2}}- \varepsilon \brac{\int_{\S^n} h_{K_{\Omega}} \, dS_{L_{\Omega,\varepsilon}} } \nonumber\\
&= \brac{W_n^{\lambda',\varepsilon} \brac{\Omega}^2+\mathrm{Area} (\partial L_{\Omega,\varepsilon})^2}^{\frac{1}{2}}- \varepsilon W_n^{\lambda', \varepsilon} \brac{\Omega}. \label{eq:Wn+1Wn-1}
\end{align}
By the isoperimetric inequality in $\R^{n+1}$, we have
\begin{equation}\label{ineq:Isoperi}
\mathrm{Area} (\partial L_{\Omega,\varepsilon}) \geq \brac{n+1} |B^{n+1}|^{\frac{1}{n+1}} \Vol \brac{L_{\Omega,\varepsilon}}^{\frac{n}{n+1}} = \brac{n+1}^{\frac{1}{n+1}}|B^{n+1}|^{\frac{1}{n+1}} W_{n+1}^{\lambda',\varepsilon} \brac{\Omega}^{\frac{n}{n+1}}.
\end{equation}
Substituting this into \eqref{eq:Wn+1Wn-1}, we arrive at 
\begin{align}
W_{n+1}^{\lambda',\varepsilon} \brac{\Omega}+\varepsilon W_n^{\lambda', \varepsilon} \brac{\Omega} \geq \brac{W_{n}^{\lambda',\varepsilon} \brac{\Omega}^2+ \brac{n+1}^{\frac{2}{n+1}} |B^{n+1}|^{\frac{2}{n+1}} W_{n+1}^{\lambda',\varepsilon} \brac{\Omega}^{\frac{2n}{n+1}} }^{\frac{1}{2}}. \label{eq:Wn+1Wn}
\end{align}
Thus, $W_{n+1}^{\lambda', \varepsilon} \brac{\Omega}+ \varepsilon W_n^{\lambda',\varepsilon} \brac{\Omega} > 0$. Squaring both sides of \eqref{eq:Wn+1Wn}, we obtain
\begin{align}
0 &\leq W_{n+1}^{\lambda', \varepsilon} \brac{\Omega}^2 - \brac{n+1}^{\frac{2}{n+1}} |B^{n+1}|^{\frac{2}{n+1}} W_{n+1}^{\lambda',\varepsilon} \brac{\Omega}^{\frac{2n}{n+1}} +2 \varepsilon W_n^{\lambda',\varepsilon} \brac{\Omega} W_{n+1}^{\lambda',\varepsilon} \brac{\Omega} \nonumber\\
&\qquad + \brac{\varepsilon^2-1} W_n^{\lambda',\varepsilon} \brac{\Omega}^2. \label{ineq:Wn+1-Wn}
\end{align}

We claim that the above inequality implies that $W_{n+1}^{\lambda', \varepsilon} \brac{\Omega}$ is greater than or equal to a suitable single variable function of $W_n^{\lambda',\varepsilon} \brac{\Omega}$. Consider a function 
\begin{equation*}
\zeta_{\alpha} \brac{t} = t^2- c_0 t^{\frac{2n}{n+1}}+ 2\varepsilon \alpha t + \brac{\varepsilon^2-1} \alpha^2, \quad t \in (-\varepsilon \alpha, \infty).
\end{equation*}
where $c_0$ and $\alpha$ are positive constants.  Note that 
\begin{equation*}
\frac{1}{2}t^{-\frac{n-1}{n+1}}\zeta_{\alpha}'\brac{t} = t^{\frac{2}{n+1}}+\varepsilon \alpha t^{-\frac{n-1}{n+1}}-\frac{n}{n+1}c_0.
\end{equation*}
The right-hand side is strictly increasing on $(-\varepsilon\alpha, \infty)$, negative near $-\varepsilon \alpha$, and tends to $ \infty$ as $t \to \infty$. Hence, it has a unique zero $t_0 \in (-\varepsilon\alpha, \infty)$. This implies that the function $\zeta_{\alpha}\brac{t}$ is strictly decreasing on $(-\varepsilon\alpha, t_0]$ and strictly increasing on $[t_0, \infty)$. Since $\zeta_{\alpha}$ is negative near $-\varepsilon \alpha$ and $\zeta_{\alpha} \brac{\infty} = \infty$, we conclude that $\zeta_{\alpha} \brac{t}$ has a unique zero in $(-\varepsilon \alpha, \infty)$ and is strictly increasing beyond this zero. Taking $c_0 = \brac{n+1}^{\frac{2}{n+1}} |B^{n+1}|^{\frac{2}{n+1}}$, $\alpha= W_n^{\lambda',\varepsilon} \brac{\Omega}$, and $t= W_{n+1}^{\lambda',\varepsilon} \brac{\Omega}$, the claim follows from the above argument and \eqref{ineq:Wn+1-Wn}. 

On the other hand, assume that \eqref{eq:Wn+1Wn} holds equality. Then \eqref{eq:Wn+1Wn-1} and \eqref{ineq:Isoperi} hold equality.
Due to the equality case of \eqref{ineq:Isoperi}, $L_{\Omega,\varepsilon}$ is a Euclidean ball, which implies that $dS_{L_{\Omega,\varepsilon}}$ is a constant multiple of the standard area measure on $\S^n$. Then, equality in \eqref{eq:Wn+1Wn-1} requires that $h_{K_{\Omega}}$ is constant, which implies that $\Omega$ is a centered geodesic ball in $\H^{n+1}$. Conversely, \eqref{eq:Wn+1Wn} holds equality in this case. Hence, we conclude
\begin{equation}\label{eq:Wn+1Wn-final}
W_{n+1}^{\lambda',\varepsilon} \brac{\Omega} \geq h_{n+1,\varepsilon} \circ h_{n, \varepsilon}^{-1} (W_n^{\lambda',\varepsilon} \brac{\Omega}).
\end{equation}
Equality holds in \eqref{eq:Wn+1Wn-final} if and only if $\Omega$ is a centered geodesic ball in $\H^{n+1}$.

\textbf{Step 2.} Next, we prove \eqref{MainIneq:ShifHuLi-H} for the cases $k=n+1$, $0\leq l \leq n$ by induction on $l$. Assume that 
\begin{equation*}
W_{n+1}^{\lambda',\varepsilon} \brac{\Omega} \geq h_{n+1,\varepsilon} \circ h_{l_0,\varepsilon}^{-1} (W_{l_0}^{\lambda',\varepsilon} \brac{\Omega})
\end{equation*}
holds for some $1 \leq l_0 \leq n$, and equality holds if and only if $\Omega$ is a centered geodesic ball in $\H^{n+1}$. The inequality \eqref{ineq:GeneralV-AFmain} implies
\begin{equation*}
W_{l_0}^{\lambda',\varepsilon}\brac{\Omega}^{\frac{n+2-l_0}{n+1-l_0}} \geq W_{l_0 -1}^{\lambda',\varepsilon} \brac{\Omega} W_{n+1}^{\lambda',\varepsilon} \brac{\Omega}^{\frac{1}{n+1-l_0}}.
\end{equation*}
Thus,
\begin{align*}
W_{l_0-1}^{\lambda',\varepsilon} \brac{\Omega} &\leq W_{n+1}^{\lambda',\varepsilon} \brac{\Omega}^{-\frac{1}{n+1-l_0}} \brac{h_{l_0,\varepsilon} \circ h_{n+1,\varepsilon}^{-1} (W_{n+1}^{\lambda',\varepsilon} \brac{\Omega} ) }^{\frac{n+2-l_0}{n+1-l_0}}\\
&= h_{l_0-1,\varepsilon} \circ h_{n+1,\varepsilon}^{-1} (W_{n+1}^{\lambda',\varepsilon} \brac{\Omega}),
\end{align*}
where we used the identity 
\begin{equation*}
h_{n+1,\varepsilon}\brac{t}^{\frac{-1}{n+1-l_0}} h_{l_0,\varepsilon} \brac{t}^{\frac{n+2-l_0}{n+1-l_0}} = h_{l_0-1,\varepsilon}\brac{t}.
\end{equation*}
The above inequality is equivalent to
\begin{equation*}
W_{n+1}^{\lambda',\varepsilon} \brac{\Omega} \geq h_{n+1,\varepsilon} \circ h_{l_0-1,\varepsilon}^{-1} (W_{l_0-1}^{\lambda',\varepsilon} \brac{\Omega}).
\end{equation*}
Consequently, we obtain \eqref{MainIneq:ShifHuLi-H} for $k=n+1$ and any $0\leq l\leq n$:
\begin{equation}\label{ineq:HLCon-k=n+1}
W_{n+1}^{\lambda',\varepsilon} \brac{\Omega} \geq h_{n+1,\varepsilon} \circ h_{l,\varepsilon}^{-1} (W_l^{\lambda',\varepsilon} \brac{\Omega}).
\end{equation}
Equality holds if and only if $\Omega$ is a centered geodesic ball in $\H^{n+1}$.

\textbf{Step 3.} Finally, for any $0 \leq l <k \leq n$, it follows from \eqref{MainIneq:AlexFen} and \eqref{ineq:HLCon-k=n+1} that
\begin{align*}
W_k^{\lambda',\varepsilon} \brac{\Omega} &\geq W_l^{\lambda',\varepsilon} \brac{\Omega}^{\frac{n+1-k}{n+1-l}} W_{n+1}^{\lambda',\varepsilon} \brac{\Omega}^{\frac{k-l}{n+1-l} } \\
&\geq W_l^{\lambda',\varepsilon} \brac{\Omega}^{\frac{n+1-k}{n+1-l}} \brac{ h_{n+1,\varepsilon} \circ h_{l,\varepsilon}^{-1} (W_l^{\lambda',\varepsilon} \brac{\Omega})}^{\frac{k-l}{n+1-l} },
\end{align*}
equivalently,
\begin{equation*}
W_k^{\lambda',\varepsilon} \brac{\Omega} \geq h_{k,\varepsilon} \circ h_{l,\varepsilon}^{-1} (W_l^{\lambda',\varepsilon} \brac{\Omega}).
\end{equation*}
Equality holds if and only if $\Omega$ is a geodesic ball in $\H^{n+1}$ centered at the origin $\mathcal{O}$. This completes the proof of Theorem \ref{Thm:FullRangeHuLi}.

\appendix 
\section{Weighted isoperimetric problem on the sphere} \label{App:A}
In this section, we show that among all domains in $\S^{n+1}_+$ with a fixed weighted volume $W_0^{\lambda'}$, the corresponding geodesic ball centered at $\mathcal{O}$ does not minimize the weighted area $W_1^{\lambda'}$. For this reason, a spherical version of Theorem \ref{Thm:FullRangeHuLi} does not hold.

To show this assertion, we first provide a formula for expressing the area element with weight $\lambda'\brac{r}$  of a hypersurface in $\S^{n+1}_+$ via the projection $\pi$; see \cite{LX23} for the hyperbolic case.
\begin{lem}
Let $\Sigma$ be a smooth hypersurface in $\S^{n+1}_+$, and let $\widehat{\Sigma} = \pi \brac{\Sigma}$. Then
\begin{equation}\label{rel-Area-spher}
\pi_* \brac{\lambda'\brac{r} \, d\mu_{\Sigma}} = \sqrt{1-u_{\widehat{\Sigma}}^2} \, d\mu_{\widehat{\Sigma}},
\end{equation}
where $u_{\widehat{\Sigma}}$ is the support function of $\widehat{\Sigma}$.
\end{lem}

\begin{proof}
Denote by $\vol_{\S}$ and $\vol_{\R}$ the volume forms of $\S^{n+1}_+$ and $\R^{n+1}$, respectively. Using \eqref{eq:rel-VolumeMeasure} and \eqref{GibS-nu}, we obtain
\begin{align*}
\pi_* \brac{\lambda'\brac{r} \, d\mu_{\Sigma}} = \pi_{*} \brac{\nu_{\Sigma} \lrcorner \brac{\lambda'\brac{r} \, d\vol_{\S}} }= \pi_{*} \nu_{\Sigma} \lrcorner \, d\vol_{\R} \\
= \sscalar{\pi_{*} \nu_{\Sigma}, \nu_{\widehat{\Sigma}}} \, d\mu_{\widehat{\Sigma}} = \sqrt{1-u_{\widehat{\Sigma}}^2} \, d\mu_{\widehat{\Sigma}},
\end{align*}
where $\lrcorner$ denotes the interior product.
\end{proof}

\begin{prop}
Let $\rho \in \brac{0, \pi/2}$.
Let $\Omega \subset \S^{n+1}_+$ be a geodesic ball of radius $\rho$ centered at $\mathcal{O}$. Then there exists a convex, static convex domain $\bar{\Omega}$ in $\S^{n+1}_+$ such that $W_0^{\lambda'}\brac{\bar{\Omega}} = W_0^{\lambda'} \brac{\Omega}$ and $W_1^{\lambda'} \brac{\bar{\Omega}} < W_1^{\lambda'}\brac{\Omega}$.
\end{prop}

\begin{proof}
Let $\theta$ be a unit vector in $\R^{n+1}$. 
Define $\Omega_t = \pi^{-1} \brac{K_{\Omega} + t \theta }$ for small $t \geq 0$, where $K_{\Omega} = \pi \brac{\Omega}$ is a centered Euclidean ball of radius $\lambda \brac{\rho}$. Then $K_{\Omega_t} = \pi \brac{\Omega_t}$ is a Euclidean ball of radius $\lambda \brac{\rho}$ centered at $t \theta$. By Proposition \ref{prop:GibS-StaticC-Convex}, $\Omega_t$ is static convex. Since $\Omega$ is strictly convex, we have $\Omega_t$ convex for small $t>0$.

It follows from \eqref{eq:Vol-WeightedVol} that $W_0^{\lambda'} \brac{\Omega_t} = W_0^{\lambda'} \brac{\Omega}$. Note that $u_{\partial K_{\Omega_t}} \brac{p+ t \theta} = \lambda \brac{\rho} + t\sscalar{\theta, \nu_{\partial K_{\Omega}} \brac{p}}$ holds for any $p \in \partial K_{\Omega}$. Then by \eqref{rel-Area-spher},
\begin{align*}
W_1^{\lambda'}\brac{\Omega_t} = \int_{\partial K_{\Omega_t}} \sqrt{1-u_{\partial K_{\Omega_t}}^2} \, d\mu_{\partial K_{\Omega_t}} = \int_{\partial K_{\Omega}} \sqrt{1-\brac{\lambda\brac{\rho} + t \sscalar{\theta, \nu_{\partial K_{\Omega}}}}^2 } \, d\mu_{\partial K_{\Omega}}.
\end{align*}
Differentiating both sides with respect to $t$ and using the divergence theorem, we have
\begin{align*}
\left. \frac{d}{dt} \right|_{t=0} W_1^{\lambda'} \brac{\Omega_t} = -\lambda'\brac{\rho}^{-1} \lambda \brac{\rho}\int_{\partial K_{\Omega}} \sscalar{\theta, \nu_{\partial K_{\Omega}}}  \, d\mu_{\partial K_{\Omega}} =0. 
\end{align*}
Moreover,
\begin{align*}
\left. \frac{d^2}{dt^2} \right|_{t=0} W_1^{\lambda'} \brac{\Omega_t} = -\lambda'\brac{\rho}^{-3} \int_{\partial K_{\Omega}} \sscalar{\theta, \nu_{\partial K_{\Omega}}}^2 \, d\mu_{\partial K_{\Omega}} <0.
\end{align*}
Consequently, we have $W_1^{\lambda'} \brac{\Omega_t} < W_1^{\lambda'} \brac{\Omega}$ for small $t>0$. The existence of the desired $\bar{\Omega}$ then follows.
\end{proof}

\section{Alternative form of the Steiner formulas in space forms}\label{App:B}
We can deduce from \eqref{eq:pullbackVolForm} that
\begin{align*}
\abs{\Sigma_t} &= \int_{\Sigma} \det \brac{\lambda\brac{t} \II_{\Sigma} + \lambda'\brac{t} g_{\Sigma}} \, d\mu_{\Sigma}\\
&= \sum_{k=0}^{n} \binom{n}{k} \lambda\brac{t}^k \lambda'\brac{t}^{n-k} \int_{\Sigma} p_k \brac{\kappa} \, d\mu_{\Sigma}.
\end{align*}
By integrating the above formula, we obtain the Steiner formula for domains in space forms with smooth boundaries \cite{San04}*{IV~18.4}:
\begin{equation}\label{eq:Steiner-SF}
\Vol \brac{\Omega_t} = \Vol \brac{\Omega} + \sum_{k=0}^n \binom{n}{k} \brac{\int_0^t \lambda \brac{s}^k \lambda'\brac{s}^{n-k} \, ds} \int_{\Sigma} p_k \brac{\kappa} d\mu.
\end{equation}
Here, $t \in [0, \delta)$ for some sufficiently small $\delta>0$.

We now recall the quermassintegrals for domains in space forms. Let $\Omega$ be a bounded domain in $\M^{n+1} \brac{c}$ with smooth boundary $\Sigma$, where $c \in \set{-1,0,1}$. The quermassintegrals $W_k$ of $\Omega$ are defined by
\begin{equation}\label{def:quermassSF}
\begin{aligned}
W_0 \brac{\Omega} &=\brac{n+1} \Vol \brac{\Omega}, \quad W_1 \brac{\Omega} = \abs{\Sigma}, \\
W_{k+1} \brac{\Omega} &= \int_{\Sigma} p_k \brac{\kappa} \, d\mu + c \frac{k}{n+2-k} W_{k-1} \brac{\Omega}, \quad k=1, \ldots, n;
\end{aligned}
\end{equation}
see \cite{Sol06}*{Prop.~7}.
This reduces to \eqref{eq:quermass-smooth} for bounded convex domains in $\R^{n+1}$. Note that $W_{n+1} \brac{\Omega} = \brac{n+1} |B^{n+1}|$ if $\Omega$ has Euler characteristic $1$ \cite{Sol06}. Under any variation $\partial_t X = \mathcal{F} \nu$, we have
\begin{equation*}
\frac{d}{dt} W_k \brac{\Omega_t} = \brac{n+1-k}\int_{\Sigma_t} \mathcal{F} p_k \brac{\kappa} \, d\mu_{\Sigma_t}, \quad k=0, \ldots,n;
\end{equation*}
see, e.g., \cite{WX14}*{Prop.~3.1}.

In view of the Steiner type formula in Proposition \ref{prop:SteinerFormulaWeightedCurvature}, it is natural to ask whether the above quermassintegrals appear as coefficients in an alternative form of the Steiner formula \eqref{eq:Steiner-SF}. We address this in the following proposition, which might be useful in future applications.
\begin{prop}
Let $\Omega$ be a bounded domain in $\M^{n+1}$ with smooth boundary. 
The Steiner formula \eqref{eq:Steiner-SF} can be expressed alternatively as
\begin{equation*}
W_0 \brac{\Omega_t} = \sum_{k=0}^n \binom{n+1}{k} \lambda'\brac{t}^{n-k} \lambda\brac{t}^k W_k \brac{\Omega} + \brac{n+1} \brac{\int_0^t \lambda\brac{s}^n \, ds } W_{n+1} \brac{\Omega}.
\end{equation*}
\end{prop}

\begin{proof}
Multiplying both sides of  \eqref{eq:Steiner-SF} by  $n+1$ and then using \eqref{def:quermassSF}, we obtain
\begin{align}
W_0 \brac{\Omega_t} &=  \sum_{k=1}^n \brac{n+1} \binom{n}{k} \brac{\int_0^t \lambda\brac{s}^k \lambda'\brac{s}^{n-k} \, ds} \brac{W_{k+1} \brac{\Omega} - c \frac{k}{n+2-k} W_{k-1} \brac{\Omega}} \nonumber\\
&\qquad + W_0 \brac{\Omega} + \brac{n+1} \brac{\int_0^t \lambda'\brac{s}^n \, ds} W_1 \brac{\Omega} \nonumber\\
&= \sum_{l=1}^{n-1} \brac{n+1} W_l\brac{\Omega} \bigg(\binom{n}{l-1} \int_0^t \lambda\brac{s}^{l-1}\lambda'\brac{s}^{n-l+1} \, ds  \nonumber\\
&\qquad - c \binom{n}{l+1} \frac{l+1}{n+1-l}\int_0^t \lambda\brac{s}^{l+1} \lambda'\brac{s}^{n-l-1} \, ds  \bigg) \nonumber\\
&\qquad + W_0 \brac{\Omega} \brac{1-cn \int_0^t\lambda\brac{s} \lambda'\brac{s}^{n-1} \, ds}+ n\brac{n+1} W_n\brac{\Omega} \int_0^t \lambda\brac{s}^{n-1} \lambda'\brac{s} \, ds \nonumber\\
&\qquad +\brac{n+1}W_{n+1} \brac{\Omega} \int_0^t \lambda\brac{s}^n \, ds. \label{eq:W0Omegat-1}
\end{align}
For $l=1, \ldots, n-1$, we note that
\begin{equation*}
\brac{n+1} \binom{n}{l-1} = \binom{n+1}{l} l, \quad 
\brac{n+1} \binom{n}{l+1} \frac{l+1}{n+1-l} = \binom{n+1}{l} \brac{n-l}.
\end{equation*}
Moreover, using integration by parts and the identity $\lambda''+ c \lambda =0$ for $\M^{n+1} = \M^{n+1} \brac{c}$, we have
\begin{equation*}
l \int_0^t \lambda\brac{s}^{l-1} \lambda'\brac{s}^{n-l+1} \, ds = \lambda\brac{t}^l \lambda'\brac{t}^{n-l} + c \brac{n-l} \int_0^t \lambda \brac{s}^{l+1} \lambda'\brac{s}^{n-l-1}  \, ds
\end{equation*}
for $l=1, \ldots, n$. Since $\lambda'\brac{0}=1$, we have
\begin{equation*}
1- cn \int_0^t \lambda\brac{s}\lambda'\brac{s}^{n-1} \, ds = \lambda'\brac{t}^n.
\end{equation*}

Substituting these identities into \eqref{eq:W0Omegat-1}, we arrive at
\begin{align*}
W_0 \brac{\Omega_t} &= \sum_{l=1}^{n} \binom{n+1}{l} W_l \brac{\Omega} \brac{l \int_0^t \lambda\brac{s}^{l-1} \lambda'\brac{s}^{n-l+1} \, ds -c \brac{n-l} \int_0^t \lambda\brac{s}^{l+1} \lambda'\brac{s}^{n-l-1} \, ds} \\
&\qquad +W_0 \brac{\Omega} \brac{1-cn \int_0^t \lambda\brac{s} \lambda'\brac{s}^{n-1} \, ds}+ \brac{n+1} W_{n+1} \brac{\Omega} \int_0^t \lambda \brac{s}^n \, ds\\
&= \sum_{l=0}^n \binom{n+1}{l} W_l \brac{\Omega} \lambda\brac{t}^{l} \lambda'\brac{t}^{n-l} + \brac{n+1} W_{n+1} \brac{\Omega} \int_0^t \lambda\brac{s}^n \, ds.
\end{align*}
This completes the proof.
\end{proof}

\end{document}